\input amstex

\documentstyle{amsppt}

\refstyle{A}

\nologo

%\magnification=\magstep1 

\hoffset .35 true in
\voffset .2 true in

\hsize=5.9 true in
\vsize=8.3 true in

\define\Adot{\bold A^\bullet}

\define\Fdot{\bold F^\bullet}
\define\Pdot{\bold P^\bullet}

\define\supp{\operatorname{supp}}

\define\gecc{\operatorname{gecc}}
\define\con{\overline{T^*_{{}_{S_\alpha}}\Cal U}}
\define\conc{\Big[ \overline{T^*_{{}_{S_\alpha}}\Cal U}\Big ]}
\define\imdf{\operatorname{im}d\tilde f}
\def\od{\overset\circ\to{\Bbb D}}
\define\dm{\operatorname{dim}}

\topmatter

\title Singularities and Enriched Cycles \endtitle

\author David B. Massey \endauthor

\address{David B. Massey, Dept. of Mathematics, Northeastern University, Boston, MA, 02115, USA} \endaddress

\email{DMASSEY\@NEU.edu}\endemail


\abstract We introduce graded, enriched characteristic cycles as a method for encoding Morse modules of strata with respect to a
constructible complex of sheaves. Using this new device, we obtain results for arbitrary complex analytic
functions on arbitrarily singular complex analytic spaces.
\endabstract

\endtopmatter

\document

\baselineskip= 12pt

\noindent\S1. {\bf Introduction}  

\vskip .1in

Let $X$ be a complex analytic space, with arbitrary singularities, and let
$\bold p\in X$. As we are interested in local questions near $\bold p$, we may assume that $X$ is a (closed) complex analytic subspace of a connected, open subset  $\Cal U\subseteq\Bbb C^{n+1}$, and
$\bold z:= (z_0, \dots, z_n)$ are coordinates on $\Bbb C^{n+1}$.  
 
Let $\tilde f:\Cal U\rightarrow\Bbb C$ be a complex analytic
function, and let $f:= {\tilde f}_{|_X}$. Let $\Cal S:=\{S_\alpha\}$ be a complex analytic Whitney stratification of $X$, with connected strata. Fix a base ring $R$ that is a regular, Noetherian ring with finite Krull dimension (e.g., $\Bbb Z, \Bbb Q, \text{or}\ \Bbb C$). Let $\Fdot$ be a bounded complex of sheaves of $R$-modules on $X$, which is constructible with respect to $\{S_\alpha\}$.

\vskip .1in

In this extremely general situation, we are interested in producing effectively calculable algebraic data which provides information about the topology and geometry of the hypersurface $V(f)$. In particular, we are interested in Thom's $a_f$ condition, and Milnor fibration data. Thus, we are led to study the nearby and vanishing cycles of $\Fdot$ along $f$; that is, respectively, $\psi_f\Fdot$ and $\phi_f\Fdot$.

\vskip .1in

Our method for producing algebraic data is motivated by our results on L\^e cycles and L\^e numbers ([{\bf M3}], [{\bf M4}], [{\bf M5}]), which corresponds to the case of an affine hypersurface where $\Fdot$ is the constant sheaf $\Bbb C^\bullet_{\Cal U}$. Much of our past work has centered around the problems of extended our work on L\^e cycles to the general setting of this paper. In [{\bf M8}], we associated cycles to arbitrary complexes of sheaves.  In [{\bf M1}], [{\bf M2}], we studied the special cases of $\psi_f\Fdot$ and $\phi_f\Fdot$. In [{\bf M7}], we gave a somewhat satisfactory generalization; we defined {\it L\^e-Vogel cycles and numbers}. 

The L\^e cycles of [{\bf M7}] were produced by starting with the characteristic cycle of the complex $\Fdot$, and then using a L\^e-Vogel process on an ideal, $J$, defining the image of $d\tilde f$ in $T^*\Cal U$. It was crucial that we had a result which told us how the characteristic cycle of $\phi_f\Fdot$ is obtained from the characteristic cycle of $\Fdot$ via blowing-up $J$; we refer to this result as the {\it vanishing index theorem}.

The fundamental weakness of the results of [{\bf M7}] lies in the fact that the characteristic cycle of $\Fdot$ contains only Euler characteristic data on the Morse modules to strata, and thus disposes of a great deal of cohomological data which describes the structure of $\Fdot$. In [{\bf M7}] and [{\bf M9}], we used the functorial properties of perverse cohomology to extract  more refined data from characteristic cycles. In [{\bf M6}], we showed with some difficulty how the microlocal theory of Kashiwara and Schapira could be used to obtain better results in the special case of an isolated critical point. 

\vskip .1in

In this paper, we show how our results from [{\bf M7}] can easily be turned into stronger results -- without using the devices of perverse cohomology or microlocal theory -- simply by using cycles with module-coefficients; these are our {\it enriched cycles}. Since we wish for our enriched cycles to be able to encode the cohomology modules of normal data to strata (with coefficients in $\Fdot$), we consider {graded, enriched cycles} and define $\gecc^\bullet(\Fdot)$, the {\it graded, enriched, characteristic cycle of $\Fdot$}. The fact that $\gecc^0\big({}^{\mu}\hskip -.02in H^k(\Fdot)\big) = \gecc^k(\Fdot)$ explains why we may use $\gecc^\bullet$ to obtain any of the results of [{\bf M7}] that used the standard characteristic cycle and perverse cohomology. 

The fundamental philosophy of this paper is: any result on perverse sheaves which is proved by using characteristic cycles and intersection theory can be generalized to arbitrary bounded, constructible complexes of sheaves by using graded, enriched characteristic cycles and a corresponding enriched intersection theory. Most of this paper consists of demonstrating applications of this general philosophy.

\vskip .1in

However, there are several important/interesting results in this paper which are not simply enriched forms of earlier results. The equivalent characterizations of {\it isolating coordinates} given in Theorem 5.10 yield a substantial improvement over [{\bf M7}] and [{\bf M8}] in the level of genericity required of our choice of coordinates when generalizing absolute polar varieties; this is very important since we are trying to produce a method for making effective computations.The result of Theorem 4.8 -- necessary and sufficient conditions for Thom's $a_f$ condition -- is important in Theorem 6.5 for similar reasons: we use 4.8 and 5.10 to produce a useful notion of ``generic coordinates'' in the relative situation. Proposition 6.11 tells us that, when the dimension of the critical locus is small, we may check whether our coordinates are generic enough ``on-the-fly''. Another interesting result which falls out of the enriched approach is a sort of converse to the theorem that: if $\Fdot$ is perverse, then $\phi_f\Fdot$ is perverse. In Corollary 3.7, we show that: if $f$ is in the square of the maximal ideal of $X$ at $\bold p$ and $\phi_f\Fdot$ is perverse in a neighborhood of $\bold p$, then $\Fdot$ is perverse in a neighborhood of $\bold p$.

\vskip .2in

The remaining sections of this paper are organized as follows.

\vskip .1in

In Section 2, we define enriched cycles, describe a corresponding enriched intersection theory (for proper intersections in smooth manifolds), and define and give some basic properties of the graded, enriched characteristic cycle of a complex of sheaves.

\vskip .1in

Section 3 contains enriched versions of a number of previous results. Of particular importance is Theorem 3.4, an enriched form of our vanishing index theorem. The only new results in this section are the corollaries to 3.4.

\vskip .1in

In Section 4, we discuss partitions, stratifications, Thom's $a_f$ condition, and {\it essentially transverse coordinates}. For most of our results, we do not need the full power of a Whitney stratification, instead we use $\Fdot$-partitions. This weakening of the type of ``stratification'' used leads to stronger results. The genericity that we use on the coordinates in Section 6 will be that the coordinates are essentially transverse to a type of partition.

\vskip .1in

In Section 5, we look at extreme generalizations of polar varieties and polar multiplicities; we define the {\it graded, enriched characteristic polar cycles} and the {\it characteristic polar modules}. We define {\it $\Fdot$-isolating coordinates} and characterize them in terms of intersections with $\gecc^\bullet(\Fdot)$. This provides the fundamental relation between enriched intersections with $\gecc^\bullet(\Fdot)$ and iterated vanishing and nearby cycles.

\vskip .1in

In Section 6, we apply the results of Section 5 to the special case of $\phi_f\Fdot$. By combining the work of Section 5 with an enriched version of the vanishing index theorem, we obtain our generalization of the L\^e cycles and L\^e numbers: the {\it graded, enriched L\^e-Vogel cycles and the L\^e-Vogel modules}. Our work on essentially transverse coordinates and the $a_f$ condition from Section 4 are used in this section to give sufficient conditions on the choice of coordinates.

\vskip .1in

In Section 7, we give a sample calculation of L\^e-Vogel cycles. We compare and contrast this with other methods of analyzing the perverse cohomology.

\vskip .1in

Section 8 contains some concluding remarks and questions.

\vskip .3in

\noindent\S2. {\bf Enriched Cycles}  

\vskip .1in

We will now define enriched cycles, graded enriched cycles, the graded enriched characteristic cycle, and operations on them.
The notions that we define are easy and/or obvious; hence, this section consists of
a series of definitions with some remarks.
 
 \vskip .1in
 
 We continue with all of the notation from the Introduction. We also need to establish some more notation that we shall use throughout the remainder of this paper.

\vskip .1in

The  cotangent bundle $T^*\Cal U @>\eta>>\Cal U$ is a trivial bundle; we will frequently use that $T^*\Cal U\cong \Cal U\times\Bbb C^{n+1}$ and write simply that $\eta$ is the projection from $\Cal U\times\Bbb C^{n+1}$ onto $\Cal U$. A linear choice of coordinates $\bold z$ for $\Cal U$ determines a basis $dz_0, \dots, dz_n$ for the cotangent spaces, and we use $\bold w:=(w_0, \dots, w_n)$ for the cotangent coordinates.

Using these coordinates, the image of $d\tilde f$ in $T^*\Cal U$, $\imdf$, is given by 
$$
V\Big(w_0-\frac{\partial\tilde f}{\partial z_0}, \dots, w_n-\frac{\partial\tilde f}{\partial z_n}\Big).
$$

As $T^*\Cal U$ is (complex) conic, we may projectivize in the cotangent directions, and consider $\Bbb P(T^*\Cal U)\cong \Cal U\times\Bbb P^n$. In this and the remaining sections of this paper, it will be important for us to consider all four projection maps: 
$$
\CD
\Cal U\times\Bbb C^{n+1}\times\Bbb P^n@>\pi>> \Cal U\times\Bbb P^n\\
@V\xi VV  @VV\nu V\\
\Cal U\times\Bbb C^{n+1} @>\eta>>\Cal U.
\endCD
$$

The conormal bundle to a stratum $S_\alpha$ is given by 
$$
T^*_{{}_{S_\alpha}}\Cal U := \big\{(\bold x, \omega)\in T^*\Cal U\ |\ \bold x\in S_\alpha, \ \omega(T_\bold x S_\alpha)\equiv 0\big\}.
$$
We are frequently interested in the closure, $\overline{T^*_{{}_{S_\alpha}}\Cal U}$, of  $T^*_{{}_{S_\alpha}}\Cal U$ inside $T^*\Cal U$. Both  $T^*_{{}_{S_\alpha}}\Cal U$ and $\overline{T^*_{{}_{S_\alpha}}\Cal U}$ are conic and may be projectivized to yield  $\Bbb P\big(T^*_{{}_{S_\alpha}}\Cal U\big)$ and $\Bbb P\big(\overline{T^*_{{}_{S_\alpha}}\Cal U}\big) = \overline{\Bbb P\big(T^*_{{}_{S_\alpha}}\Cal U\big)}$.

\vskip .1in

Finally, it is fundamental throughout this paper that the functors $\psi_f$ and $\phi_f$, shifted by $-1$, commute, up to natural isomorphism, with (middle perversity) perverse cohomology ${}^{\mu}\hskip -.02in H^*$; see [{\bf BBD}], 10.3.13 of [{\bf K-S2}] and Remark 6.0.6 of Section 6.0.4 of [{\bf S}] (but, be aware that the definition of the vanishing cycles used in [{\bf K-S2}] is shifted by $1$ from the standard definition that we use).  In particular, the shifted nearby and vanishing cycles are {\it perverse functors}, i.e., they take perverse sheaves to perverse sheaves. For the nearby cycles, this was first proved by Goresky and MacPherson in [{\bf G-M1}]. The first proof of which we are aware that the shifted vanishing cycle functor is perverse appears in 1.7 of [{\bf Br}]. For these reasons, we shall always include a shift of $-1$ when we use the nearby and vanishing cycle functors; we shall write $\psi_f[-1]$ and $\phi_f[-1]$ for the functors $\psi_f$ and $\phi_f$ composed with a shift by $-1$.

\vskip .3in

\noindent{\bf Definition 2.1}.  An {\it enriched cycle}, $E$, in $X$ is a formal, locally finite sum $\sum_V E_V[V]$, where the
$V$'s are irreducible analytic subsets of $X$ and the
$E_V$'s are finitely-generated $R$-modules. We refer to the $V$'s as the {\it components} of $E$, and to $E_V$ as the {\it
$V$-component module of
$E$}. Two enriched cycles are considered the same provided that all of the component modules are isomorphic. The underlying set
of $E$ is
$|E|:=
\cup_{{}_{E_V\neq 0}}V$.

If $C = \sum n_V[V]$ is an ordinary positive cycle in $X$, i.e., all of the
$n_v$ are non-negative integers, then there is a corresponding enriched cycle $[C]^{\operatorname{enr}}$ in which the
$V$-component module is the free
$R$-module of rank $n_V$. If $R$ is an integral domain, so that rank of an $R$-module is well-defined, then an enriched cycle
$E$ yields an ordinary cycle $[E]^{\operatorname{ord}}:=\sum_V (\operatorname{rk}(E_V))[V]$.

If $q$ is a finitely-generated module and $E$ is an enriched cycle, then we let
$qE:=\sum_V(q\otimes E_V)[V]$; thus, if $R$ is an integral domain and $E$ is an enriched cycle,
$[qE]^{\operatorname{ord}}=(\operatorname{rk}(q))[E]^{\operatorname{ord}}$ and if
$C$ is an ordinary positive cycle and $n$ is a positive integer, then $[nC]^{\operatorname{enr}}=R^n[C]^{\operatorname{enr}}$.

\vskip .1in

The (direct) sum of two enriched cycles $D$ and $E$ is given by $(D + E)_V := D_V\oplus E_V$. 

\vskip .1in

There is a partial ordering on enriched cycles given by: $D\leqslant E$ if and only if there exists an enriched cycle $P$ such that $D+P=E$. This relation is clearly reflexive and transitive; moreover, anti-symmetry follows from the fact that if $M$ and $N$ are Noetherian modules such that $M\oplus N\cong M$, then $N=0$.

\vskip .1in

If two irreducible analytic subsets $V$ and
$W$ intersect properly in $\Cal U$, then the (ordinary) intersection cycle $[V]\cdot[W]$ is a well-defined positive cycle; we
define the enriched intersection product of $[V]^{\operatorname{enr}}$ and $[W]^{\operatorname{enr}}$ by
$[V]^{\operatorname{enr}}\odot[W]^{\operatorname{enr}} = ([V]\cdot[W])^{\operatorname{enr}}$. If $D$ and $E$ are enriched
cycles, and every component of $D$ properly intersects every component of $E$ in $\Cal U$, then we say that {\it $D$ and $E$
intersect properly} in $\Cal U$ and we extend the intersection product linearly, i.e., if $D=\sum_V D_V[V]$ and $E=\sum_W
E_W[W]$, then
$$ D\odot E:= \sum_{V, W} (D_V\otimes E_W)([V]\cdot [W])^{\operatorname{enr}}.
$$

Suppose that $\bold g:=(g_0, \dots, g_{d})$ is a $(d+1)$-tuple of analytic functions on $\Cal U$ and the $D=\sum D_V[V]$ is an
enriched cycle in $\Cal U$. Then, for each $V$, the blow-up of the ideal generated by the restriction of $\bold g$ to $V$,
$\operatorname{Bl}_\bold g (V)$, yields an irreducible analytic subset of $\Cal U\times\Bbb P^d$; in addition, we obtain an
exceptional divisor $\operatorname{Ex}_\bold g(V)$, which is an (ordinary) positive cycle in $\Cal U\times\Bbb P^d$. We define
the {\it blow-up and exceptional divisor of $D$ along $\bold g$} to be the enriched cycles in $\Cal U\times\Bbb P^d$ given by 
$\operatorname{Bl}_\bold g (D):=\sum_V D_V[\operatorname{Bl}_\bold g (V)]^{\operatorname{enr}}$ and $\operatorname{Ex}_\bold g (D):=\sum_V D_V[\operatorname{Ex}_\bold g (V)]^{\operatorname{enr}}$, respectively. If we had started with an ordinary cycle $D$, we would analogously obtained ordinary cycles for $\operatorname{Bl}_\bold g (D)$ and $\operatorname{Ex}_\bold g (D)$.

\vskip .1in

A {\it graded, enriched cycle $E^\bullet$} is simply an enriched cycle $E^i$ for $i$ in some bounded set of integers. An single
enriched cycle is considered as a graded enriched cycle by being placed totally in degree zero. The analytic set $V$ is a {\it
component of $E^\bullet$} if and only if $V$ is a component of $E^i$ for some $i$, and the underlying set of $E^\bullet$ is
$|E^\bullet|=\cup_i|E^i|$. If $k$ is an integer, we define the {\it
$k$-shifted graded, enriched cycle}
$E^\bullet[k]$ by $(E^\bullet[k])^i:= E^{i+k}$.

If $R$ is a domain, then $E^\bullet$ yields an ordinary cycle
$$[E^\bullet]^{\operatorname{ord}}:=\sum_{i}(-1)^i[E^i]^{\operatorname{ord}} \ =\ \sum_{i, V}(-1)^i(\operatorname{rk}(E^i_V))[V].$$ 
If $q$ is a finitely-generated module and $E^\bullet$ is a graded enriched cycle, then we define the graded enriched cycle
$qE^\bullet$ by
$(qE^\bullet)^i:=\sum_V(q\otimes E^i_V)[V]$. The (direct) sum of two graded enriched cycles $D^\bullet$ and $E^\bullet$ is given
by $(D^\bullet+ E^\bullet)^i_V := D^i_V\oplus E^i_V$. If $D^i$ properly intersects $E^j$ for all $i$ and $j$, then we say that
$D^\bullet$ and $E^\bullet$ {\it intersect properly} and we define the intersection product by 
$$ (D^\bullet\odot E^\bullet)^k:=\sum_{i+j = k}(D^i\odot E^j).
$$

Let $\tau: W\rightarrow Y$ be a proper morphism between analytic spaces. If $C = \sum n_V[V]$ is an ordinary positive cycle in
$W$, then the proper push-forward $\tau_*(C) = \sum n_V\tau_*([V])$ is a well-defined ordinary cycle. If $E^\bullet=\sum_V
E^\bullet_V[V]$ is an enriched cycle in $W$, then we define the {\it proper push-forward of $E^\bullet$ by $\tau$} to be the
graded enriched cycle $\tau^\bullet_*(E^\bullet)$ defined by
$$
\tau^j_*(E^\bullet) \ := \ \sum_V E^j_V[\tau_*([V])]^{\operatorname{enr}}.
$$ The ordinary projection formula for divisors ([{\bf F}], 2.3.c) immediately implies the following enriched version. Let
$E^\bullet$ be a graded enriched cycle in $X$. Let $W:=|E^\bullet|$. Let $\tau:W\rightarrow Y$ be a proper morphism, and let
$g: Y\rightarrow\Bbb C$ be an analytic function such that
$g\circ\tau$ is not identically zero on any component of $E^\bullet$. Then, $g$ is not identically zero on any component of
$\tau^\bullet_*(E^\bullet)$ and
$$
\tau^\bullet_*\big(E^\bullet\ \odot \ V(g\circ\tau)\big) \ = \ \tau^\bullet_*(E^\bullet) \ \odot \ V(g).
$$

\vskip .1in

The standard intersection result on conservation of number generalizes easily to:

\vskip .1in

\noindent{\bf conservation of module}:

\vskip .1in

Let $E$ be a purely $k$-dimensional enriched cycle in $\Cal U$ and let $\bold f:=(f_1, \dots, f_k)\in\left(\Cal
O_{\Cal U}\right)^k$ be such that $E$ and $V(\bold f)$ intersect properly in the isolated point $\bold p$. 

  Let  
$g_1(\bold z, t),
\dots, g_k(\bold z,t)\in
\Cal O_{\Cal U\times\od}$ be such that $g_i(\bold z, 0)=f_i(\bold z)$ for all $i$. For $t_0\in\od$, let $C_{t_0}$ be the cycle in $\Cal U$ given by $\big[V(g_1(\bold z, t_0),
\dots, g_k(\bold z, t_0))\big]$.  Then,
$$
 (E\odot V(\bold f))_{\bold p} \cong\sum_{\bold q\in{\overset\circ\to{B_\epsilon}}\cap |E|\cap |C_{t_0}|} \hskip -.2in\big(E\odot
C_{t_0}\big)_{\bold q} ,
$$ where $\epsilon>0$ is sufficiently small, $\overset\circ\to{B_\epsilon}$ is an open ball of radius $\epsilon$
centered at $\bold p$, and
$|t_0|\ll\epsilon$.

\vskip .35in

\noindent{\bf Definition 2.2}.   If $D^\bullet$ is a graded, enriched cycle and $\bold g:=(g_0, \dots, g_{d})$ is a
$(d+1)$-tuple of analytic functions on $\Cal U$, then we define the {\it blow-up and exceptional divisor of $D^\bullet$ along
$\bold g$} to be the graded, enriched cycles in $\Cal U\times\Bbb P^d$ given by 
$(\operatorname{Bl}_\bold g (D^\bullet))^i:=\sum_V D^i_V[\operatorname{Bl}_\bold g (V)]^{\operatorname{enr}}$ and
$(\operatorname{Ex}_\bold g (D^\bullet))^i:=\sum_V D^i_V[\operatorname{Ex}_\bold g (V)]^{\operatorname{enr}}$, respectively.

\vskip .1in

Later, we will be especially interested in the case where $\bold g$ is the tuple defining $\imdf$, $$\Big(w_0-\frac{\partial\tilde f}{\partial z_0}, \dots, w_n-\frac{\partial\tilde f}{\partial z_n}\Big),$$  and $D^\bullet$ is such that $|D^\bullet|$ is a union of components of the form $\con$. In this case, we will denote the blow-up and exceptional divisor by $\operatorname{Bl}_{\imdf} (D^\bullet)$ and
$\operatorname{Ex}_{\imdf} (D^\bullet)$, respectively.

\vskip .3in

\noindent{\bf Definition 2.3}.  Suppose that $\Fdot$ is a bounded complex of sheaves, which is constructible with respect to an
analytic Whitney stratification $\{S_\alpha\}$, in which the strata are connected. Let $d_\alpha:=\dim S_\alpha$. If $(\Bbb
N_\alpha, \Bbb L_\alpha)$ is a pair consisting of a normal slice and complex link (see [{\bf G-M2}]), respectively, to the stratum $S_\alpha$,
then, for each integer $k$, the isomorphism-type of the module
$\Bbb H^{k-d_\alpha}(\Bbb N_\alpha, \Bbb L_\alpha; \Fdot)$ is independent of the choice of $(\Bbb N_\alpha, \Bbb L_\alpha)$; we
refer to
$\Bbb H^{k-d_\alpha}(\Bbb N_\alpha,\Bbb L_\alpha; \Fdot)$ as the {\it degree $k$ Morse module of $S_\alpha$ with respect to
$\Fdot$}. The {\it graded, enriched characteristic cycle of $\Fdot$ in the cotangent bundle $T^*\Cal U$} is defined in degree
$k$ to be 
$$
\gecc^k(\Fdot) :=\sum_\alpha H^{k-d_\alpha}(\Bbb N_\alpha,\Bbb L_\alpha; \Fdot)\big[\,\overline{T^*_{{}_{S_\alpha}}\Cal
U}\,\big].
$$

\vskip .2in

\noindent{\it Remark 2.4}. There are no canonical choices for defining the the normal slices or complex links of strata (again, see [{\bf G-M2}]). However, as two enriched cycles are equal provided that the component modules are all isomorphic, the graded, enriched characteristic cycle is well-defined.

Note that the ordinary characteristic cycle, $\operatorname{Ch}(\Fdot)$, of $\Fdot$ is related to the
graded, enriched characteristic cycle of $\Fdot$ by $\operatorname{Ch}(\Fdot)= (-1)^{\dim
X}[\gecc^\bullet(\Fdot)]^{\operatorname{ord}}$, provided that the base ring is an integral domain.

On a different note, recall that a complex $\Fdot$ is called {\it pure}  (of shift $0$) provided that $\gecc^\bullet(\Fdot)$ is concentrated in degree $0$ (see [{\bf K-S1}], 7.2 and 9.5). If $\Pdot$ is a perverse sheaf on $X$, then
$\Pdot$ is pure; the converse of this is also true (see 9.5.2 of [{\bf K-S1}] or Remark 2.6).

Finally, note that, for all $k$,  $\gecc^k(\Fdot[i]) = (\gecc^k(\Fdot))[i]$.

\vskip .3in

\noindent{\bf Proposition 2.5}. {\it For all $k$, there is an equality of enriched cycles given by $$\gecc^0\big({}^{\mu}\hskip -.02in H^k(\Fdot)\big) = \gecc^k(\Fdot).$$
Furthermore, there are equalities of sets given by $\supp{}^{\mu}\hskip -.02in H^k(\Fdot) =\eta(| \gecc^k(\Fdot)|)$,
and $\supp(\Fdot) =\eta(| \gecc^\bullet(\Fdot)|)$.

}

\vskip .2in

\noindent{\it Proof}. The first equality follows immediately from the fact that perverse cohomology commutes with taking the shifted vanishing cycles; the argument is precisely the same as that of Proposition 2.1 in  [{\bf M9}].

If $\Pdot$ is a perverse sheaf, then $\supp(\Pdot) = \eta(|\gecc^0(\Pdot)|)$; this is essentially the last part of Lemma 3.1 of [{\bf M1}], and the proof is the same. The second equality of the proposition follows by applying this to ${}^{\mu}\hskip -.02in H^k(\Fdot)$, and then using the first equality.

As $\supp(\Fdot) =\bigcup_k \supp{}^{\mu}\hskip -.02in H^k(\Fdot)$ (see [{\bf K-S2}]), the last equality of the proposition follows from the second. 

Alternatively, the last two equalities could be concluded from the fact that, at a generic point in an irreducible component of $\supp(\Fdot)$, the Morse modules simply yield stalk cohomology.\qed

\vskip .28in

\noindent{\it  Remark 2.6}. The equality of $\supp{}^{\mu}\hskip -.02in H^k(\Fdot)$ and  $\eta(| \gecc^k(\Fdot)|)$ implies that pure complexes have ${}^{\mu}\hskip -.02in H^k(\Fdot)=0$ for $k\neq 0$, i.e., we recover the result of 9.5.2 of [{\bf K-S1}] that pure complexes are perverse.

We should also remark here that the equality of $\supp(\Fdot)$ and $\eta(| \gecc^\bullet(\Fdot)|)$ explains why a number of our later results have hypotheses involving $| \gecc^\bullet(\Fdot)|$ even though we wish to conclude results in individual degrees.

\vskip .28in

\noindent{\bf Definition 2.7}.  An irreducible subvariety $Y$ of $X$ is an {\it essential subvariety for $\Fdot$}, or an {\it $\Fdot$-essential subvariety}, provided that there is a irreducible component $C$ of $|\gecc^\bullet(\Fdot)|$ such that $Y=\eta(C)$.

A connected submanifold $M\subseteq\Cal U$ is  an {\it essential submanifold for $\Fdot$}, or an {\it $\Fdot$-essential submanifold}, provided that there is an $\Fdot$-essential subvariety $Y$ such that $M=Y_{\operatorname{reg}}$.

A stratum $S_\alpha$ is {\it $\Fdot$-visible} provided that $\overline{S_\alpha}$ is an $\Fdot$-essential subvariety. This is equivalent to saying that $S_\alpha$ has a non-zero Morse module,
with respect to $\Fdot$, in some degree, i.e., provided that $\overline{T^*_{{}_{S_\alpha}}\Cal U}$ is a component of
$\gecc^\bullet(\Fdot)$ or, equivalently, that $\overline{T^*_{{}_{S_\alpha}}\Cal U}$ is an irreducible component of
$|\gecc^\bullet(\Fdot)|$.

\vskip .28in

\noindent{\it Remark 2.8}.  In terms of the definitions in 2.7, the last equality of Proposition 2.5 can be restated as $$\supp(\Fdot)\ =\bigcup\Sb\Fdot-\text{essential}\\\text{subvarieties } Y \endSb \hskip -.2in Y\ =\bigcup\Sb\Fdot-\text{essential}\\\text{submanifolds } M\endSb \hskip -.2in M\ =\bigcup_{\Fdot-\text{visible } S_\alpha}\hskip -.2in \overline{S_\alpha}.$$

\vskip .3in

\noindent\S3. {\bf Enriched Forms of Previous Results}

\vskip .2in

In this section, we describe the enriched versions of a number of known results. We provide no proofs, since all one has to do
is rewrite previous Morse-theoretic proofs, using enriched cycles in place of ordinary cycles.

We continue with all of our notation from the previous two sections.

\vskip .3in

\noindent{\bf Definition 3.1}.  If $M$ is an analytic submanifold of $\Cal U$ and $M\subseteq X$, then the
{\it  relative conormal space (of $M$ with respect to $f$ in $\Cal U$)}, $T^*_{f_{|_M}}\Cal U$, is given by
$$ T^*_{f_{|_M}}\Cal U:=\{(\bold x, \omega)\in T^*\Cal U\ |\ \bold x\in M,\ \omega\big(\operatorname{ker}d_\bold
x(f_{|_M})\big)=0\} = 
$$
$$
\{(\bold x, \omega)\in T^*\Cal U\ |\ \bold x\in M,\ \omega\big(T_\bold x M\cap\operatorname{ker}d_\bold x\tilde f\big)=0\}.
$$

If $\gecc^k(\Fdot)=\sum E^k_\alpha\,\conc$, then we define the {\it relative graded enriched conormal cycle,
$\big(T^*_{{}_{f,
\Fdot}}\Cal U\big)^\bullet$, of $f$ with respect to $\Fdot$} by 
$$\big(T^*_{{}_{f, \Fdot}}\Cal U\big)^k :=
\sum_{f_{|_{S_\alpha}}\neq\text{constant}}E^k_\alpha \Big [\overline{T^*_{f_{|_{S_\alpha}}}\Cal U}\Big ].$$

\vskip .1in

The following result is an enriched version of Theorem 2.3 of [{\bf M1}]; it follows trivially by ``enriching'' the cycles in the
Morse theoretic proof given in [{\bf M1}]. Alternatively, one can use the device of perverse cohomology, as in [{\bf M9}], to obtain the result directly from the statement of Theorem 2.3 of [{\bf M1}].

\vskip .2in

\noindent{\bf Theorem 3.2}. {\it There is an equality of graded enriched cycles given by
$$
\gecc^k(\psi_f[-1]\Fdot) \ = \ \big(T^*_{{}_{f, \Fdot}}\Cal U\big)^k\odot(V(f)\times\Bbb C^{n+1}).
$$ 

In particular, $\dsize\supp(\psi_f[-1]\Fdot) \ = \ \Big(V(f)\cap\bigcup\Sb\Fdot-\text{essential }\\\text{subvarieties }Y\\Y\not\subseteq V(f)\endSb Y\Big) 
\ = \ V(f)\cap\big(\overline{\supp(\Fdot)-V(f)}\big)$.
}

\vskip .3in

The next theorem is the main result of [{\bf M6}], stated in a form that uses our current terminology. One could also obtain this result by using enriched cycles throughout the proof in [{\bf L}].

\vskip .2in

\noindent{\bf Theorem 3.3}. {\it Suppose that $f(\bold p)=v$. 

Then,  $\dim_\bold p (\supp\phi_{f-v}[-1]\Fdot)\leqslant 0$ if and only if
$\dim_{\bold p}\eta\big(|\gecc^\bullet(\Fdot)|\cap \operatorname{im}d\tilde f\big)\leqslant 0$, and when this
is the case, $\dim_{(\bold p, d_\bold p\tilde f)}\big(|\gecc^\bullet(\Fdot)|\cap \operatorname{im}d\tilde f\big)\leqslant 0$ and 
$$ H^k(\phi_{f-v}[-1]\Fdot)_\bold p\cong \big(\gecc^k(\Fdot)\ \odot\ \operatorname{im}d\tilde f\big)_{(\bold p, d_\bold p\tilde
f)},
$$ where $\operatorname{im}d\tilde f$ is considered as a graded, enriched cycle. }

\vskip .3in

The next theorem follows from
Proposition 2.2 and Theorem 3.1 of [{\bf M9}], but we could have obtained the result directly by using enriched cycles throughout the proof of Theorem 2.10 of [{\bf M1}].

\vskip .2in

\noindent{\bf Theorem 3.4}. {\it There is an equality of closed subsets of $X$ given by
$$\bigcup_{v\in\Bbb C}\supp\phi_{f-v}[-1]\Fdot\ =\ \eta\big(|\gecc^\bullet(\Fdot)|\cap\operatorname{im}d\tilde f\big),
$$ 
and, for all $k$, an equality of graded, enriched cycles given by
$$
\sum_{v\in\Bbb C}\Bbb
P\big(\gecc^k(\phi_{f-v}[-1]\Fdot)\big)\ =\ \pi_*\big({\operatorname{Ex}}_{\operatorname{im}d\tilde f}(\gecc^k(\Fdot))\big).
$$ 
In particular, for all $k$, there is an equality of sets
$$\bigcup_{v\in\Bbb C}\eta\big(|\gecc^k(\phi_{f-v}[-1]\Fdot)|\big)\ =\ \eta\big(|\gecc^k(\Fdot)|\cap\operatorname{im}d\tilde f\big).
$$ 
}

\vskip .1in

Theorem 3.3 can be recovered very quickly as a special case of Theorem 3.4.

\vskip .2in

The following corollary is immediate from Theorem 3.4.

\vskip .2in

\noindent{\bf Corollary 3.5}. {\it Suppose that $f(\bold p)=0$, and that $M$ is an $R$-module.

 If  the enriched cycle $M\big[\{\bold p\}\times\Bbb C^{n+1}\big]$ is a summand of $\gecc^k(\Fdot)$, then $M\big[\{\bold p\}\times\Bbb C^{n+1}\big]$ is a summand of $\gecc^k(\phi_f[-1]\Fdot)$. In particular, if $\{\bold p\}\times\Bbb C^{n+1}$ is not a component of $\big|\gecc^\bullet(\phi_f[-1]\Fdot)\big|$, then it is not a component of $\big|\gecc^\bullet(\Fdot)\big|$.
}

\vskip .3in

\noindent{\bf Corollary 3.6}. {\it Suppose that $\phi_f[-1]\Fdot = 0$. 

\vskip .1in

Then, $\psi_f[-1]\Fdot\cong \Fdot_{|_{V(f)}}[-1]$ and
$$\supp(\psi_f[-1]\Fdot) = \supp( \Fdot_{|_{V(f)}}[-1])=V(f)\cap\supp \Fdot;$$
in addition, $\dim\big(V(f)\cap\supp \Fdot\big)\leqslant -1+\dim(\supp\Fdot)$. (In particular, if $\dim(\supp\Fdot)=0$, then $V(f)\cap \supp\Fdot=\emptyset$.)
}

\vskip .2in

\noindent{\it Proof}. The isomorphism is immediate from the fundamental distinguished triangle relating the nearby and vanishing cycles. By 3.4, since $\phi_f[-1]\Fdot = 0$, $V(f)$ can not contain an $\Fdot$-visible stratum. Therefore, by the final statement of 3.2, $\supp(\psi_f[-1]\Fdot)= V(f)\cap\supp \Fdot$, and the dimension of this must drop.\qed

\vskip .3in

\noindent{\bf Corollary 3.7}. {\it Suppose that $f$ is in the square of the maximal ideal of $X$ at $\bold p$, i.e., $f\in{\frak m}^2_{{}_{X, \bold p}}$. If $\phi_f[-1]\Fdot$ is perverse in a neighborhood of $\bold p$, then $\Fdot$ is perverse in a neighborhood of $\bold p$. If $\bold p$ is not in the support of $\phi_f[-1]\Fdot$, then $\bold p$ is not in the support of $\Fdot$.
}

\vskip .2in

\noindent{\it Proof}. Since $f\in{\frak m}^2_{{}_{X, \bold p}}$, it follows that, for all $S_\alpha$, if $\bold p\in\overline{S_\alpha}$, then $(\bold p, d_\bold p\tilde f)\in\con$. Therefore, Theorem 3.4 implies that, if $\gecc^\bullet(\phi_f[-1]\Fdot)$ is concentrated in degree $0$, then $\gecc^\bullet(\Fdot)$ must be concentrated in degree $0$ over a neighborhood of $\bold p$. This yields the first statement. The second statement uses the same argument.
\qed

\vskip .3in

As our final result of this section, we need to state the enriched version of Theorem I.2.20 of [{\bf M7}]. We will use this result in Section 6, where it will enable us to actually perform calculations. We state the result in the form in which we shall use it. Recall the definition of the graded, enriched exceptional divisor from 2.2.

\vskip .3in

\noindent{\bf Theorem 3.8}. {\it Let $C^\bullet$ be a purely $(n+1)$-dimensional graded, enriched cycle on $\Cal U\times\Bbb C^{n+1}$. Let $\bold h:=(h_0, \dots, h_n)$ be an $(n+1)$-tuple of analytic functions on $\Cal U\times\Bbb C^{n+1}$.

Suppose that, for all $j$ such that $0\leqslant j\leqslant n+1$, ${\operatorname{Ex}}^\bullet_{\bold h}(C^\bullet)$ properly intersects $\Cal U\times\Bbb C^{n+1}\times\Bbb P^j\times\{\bold 0\}$ inside $\Cal U\times\Bbb C^{n+1}\times\Bbb P^n$, and define ${}^k\Delta^j_\bold h$ by
$$
{}^k\Delta^j_\bold h=\xi_*\big({\operatorname{Ex}}^k_{\bold h}(C^\bullet) \ \odot \ (\Cal U\times\Bbb C^{n+1}\times\Bbb P^j\times\{\bold 0\})\big).
$$

\vskip .1in

Then, there exists an open neighborhood of $V(\bold h)$ in $\Cal U\times\Bbb C^{n+1}$ in which ${}^k\Delta^j_\bold h$ can be calculated via the following inductive process, in which all of the intersections are proper:

\vskip .1in

If $C^k =\sum_V C_V[V]$, then let $\dsize{}^k\Pi^{n+1}_\bold h:=\sum_{V\not\subseteq V(\bold h)}C_V[V]$. Then, one can write $${}^k\Pi^{n+1}_\bold h\odot V(h_n) = {}^k\Pi^{n}_\bold h+{}^k\Delta^{n}_\bold h,$$
where ${}^k\Pi^{n}_\bold h$ denotes the sum of all those components (with their component modules) of the intersection which are not contained (as sets) in $V(\bold h)$, and ${}^k\Delta^{n}_\bold h$ will precisely equal the sum of the components of the intersection which are contained in $V(\bold h)$. Proceeding inductively, one writes the intersection
$$
{}^k\Pi^{j+1}_\bold h\odot V(h_j) = {}^k\Pi^{j}_\bold h+{}^k\Delta^{j}_\bold h,
$$
where ${}^k\Pi^{j}_\bold h$ denotes the sum of all those components of the intersection which are not contained in $V(\bold h)$, and ${}^k\Delta^{i}_\bold h$ will precisely equal the sum of the components of the intersection which are contained in $V(\bold h)$. 
}

\vskip .3in

\noindent\S4. {\bf Partitions and Stratifications}  

\vskip .1in

In defining $\gecc^\bullet(\Fdot)$, we used a Whitney stratification of $X$ with respect to which $\Fdot$ was constructible. However, in practice, we do not explicitly need Whitney's condition b) or the condition of the frontier. Since weaker hypotheses on the ``stratifications'' will yield stronger, more useful results, we wish to describe the types of partitions of $X$ that we will use, and prove a few basic results.

\vskip .1in

We continue with all of the notation from the previous sections {\bf except that $\{S_\alpha\}$ is no longer assumed to be a Whitney stratification}. Throughout the remainder of this paper, if $Y$ is an analytic set, then when we write $\dim_\bold p Y\leqslant 0$, we mean that either $\bold p$ is an
isolated point of
$Y$ ($\dim_\bold p Y = 0$) or $\bold p\not\in Y$ ($\dim_\bold p Y=-\infty$).

\vskip .3in

\noindent{\bf Definition 4.1}. A {\it (complex analytic) partition of $X$} is a locally-finite collection $\Cal S:=\{S_\alpha\}_{\alpha}$ of disjoint, connected, complex analytic manifolds whose union is all of $X$ and such that, for all $\alpha$, $\overline{S_\alpha}$ and $\overline{S_\alpha}-S_\alpha$ are complex analytic subsets of $\Cal U$. We still refer to the elements of $\Cal S$ as {\it strata}.

\vskip .1in

The strata of a partition are partially-ordered by: $S_\beta\leqslant S_\alpha$ if and only if $S_\beta\subseteq\overline{S_\alpha}$.

\vskip .1in 

A partition $\Cal S$ of $X$ is a {\it stratification} provided that the condition of the frontier holds, i.e., for all $S_\alpha\in\Cal S$, $\overline{S_\alpha}$ is a union of strata. This means that $S_\beta\leqslant S_\alpha$ if and only if $S_\beta\cap\overline{S_\alpha}\neq\emptyset$.

\vskip .1in

A partition $\Cal S$ of $X$ is a {\it Whitney a) partition} provided that all pairs of strata satisfy Whitney's condition a); in conormal terms this means that, for all $S_\beta\in\Cal S$,
$$
\overline{T^*_{S_\beta}\Cal U}\ \subseteq\ \bigcup_{\alpha}T^*_{S_\alpha}\Cal U.
$$
In particular, this implies that $\bigcup_{\alpha}T^*_{S_\alpha}\Cal U$ is equal to $\bigcup_{\alpha}\con$ and, hence, is closed.

\vskip .1in

If $\Fdot$ is a bounded, constructible complex of sheaves on $X$, then a partition $\Cal S$ of $X$ is an {\it $\Fdot$-partition of $X$} provided that 
$$
|\gecc^\bullet(\Fdot)|\ \subseteq\ \bigcup_{\alpha}\con.
$$

In an $\Fdot$-partition of $X$, a {\it visible stratum}, $S_\alpha$, is one  such that $\con\subseteq |\gecc^\bullet(\Fdot)|$. Note that this extends our definition in 2.7 to the case where $\Cal S$ need not be a Whitney stratification. Also note that, if $\Cal S$ is an $\Fdot$-partition of $X$, then, as before, Proposition 2.5 implies that  $$\supp(\Fdot) =\bigcup_{\Fdot-\text{visible } S_\alpha}\hskip -.2in \overline{S_\alpha};$$ in particular, this implies that the maximal elements of $\{S_\alpha\in\Cal S\ |\ S_\alpha\subseteq \supp(\Fdot)\}$ are all visible, and that their union is equal to $\supp(\Fdot)$.

\vskip .3in

Certainly, a Whitney stratification of $X$, with connected strata, with respect to which $\Fdot$ is constructible, is an  $\Fdot$-partition of $X$. However, in our results, it is only the properties of $\Fdot$-partitions that we actually use.

\vskip .3in

\noindent{\bf Definition 4.2}. Suppose that $\Cal S$ is a partition of $X$. Let $S_\alpha\in\Cal S$, and let $\bold p\in X$ (but $\bold p$ is not necessarily in $S_\alpha$). 

The germ at $\bold p$ of an analytic submanifold $N\subseteq \Cal U$ is {\it essentially transverse to $S_\alpha$ at $\bold p$} (in $\Cal U$) provided that there is an open neighborhood of $\bold p$ in which $N$ transversely intersects $S_\alpha-\{\bold p\}$ inside $\Cal U$. 

We say that the germ $N$ is  {\it essentially transverse to $\Cal S$ at $\bold p$} provided that $N$ is essentially transverse to $S_\alpha$ at $\bold p$, for all $S_\alpha\in\Cal S$.

If $N$ is essentially transverse to $\Cal S$ at $\bold p$, then, in a neighborhood of $\bold p$, there is an {\it induced partition of $X\cap N$} in which the strata are the point-stratum $\{\bold p\}$, together with the connected-components  of the intersection $(S_\alpha-\{\bold p\})\cap N$.

The coordinates $\bold z$ for $\Cal U$ are {\it essentially transverse to $S_\alpha$ at $\bold p$} provided that, for all $i$ with $0\leqslant i\leqslant n$, $V(z_0-p_0, \dots, z_i-p_i)$ is essentially transverse to $S_\alpha$ at $\bold p$. This is clearly equivalent to requiring that, for all $i$ with $0\leqslant i\leqslant -1+\dm S_\alpha$, $V(z_0-p_0, \dots, z_i-p_i)$ is essentially transverse to $S_\alpha$ at $\bold p$.

The coordinates $\bold z$ for $\Cal U$ are {\it essentially transverse to $\Cal S$ at $\bold p$} provided that $\bold z$ is essentially transverse to each $S_\alpha$ at $\bold p$. This is equivalent to: $V(z_0-p_0)$ is essentially transverse to $\Cal S$ at $\bold p$ and, for all $i$ such that $1\leqslant i\leqslant n$, $V(z_i-p_i)$ is essentially transverse to the induced partition of $X\cap V(z_0-p_0, \dots, z_{i-1}-p_{i-1})$.

\vskip .3in

In the proposition below, we should point out that the only difference between the statements in b) and c) is that, in b), $V(z_i-p_i)$ explicitly appears in the intersection, while it does not in c).

\vskip .3in

\noindent{\bf Proposition 4.3}. {\it Suppose that $\Cal S$ is a partition of $X$. Let $S_\alpha\in\Cal S$, and let $\bold p\in X$. Then, the following are equivalent:

\vskip .1in

\noindent a) The coordinates $\bold z$ are essentially transverse to $S_\alpha$ at $\bold p$;

\vskip .1in

\noindent b) there exists an open neighborhood $\Cal W$ of $\bold p$ in $\Cal U$ such that, for all $i$ such that $0\leqslant i\leqslant n$,
$$
\Bbb P(T^*_{S_\alpha}\Cal U)\ \cap \ \Big(\big(\Cal W\cap V(z_0-p_0, \dots, z_i-p_i)\big)\times\Bbb P^i\times\{\bold 0\}\Big) \ \subseteq \ \{\bold p\}\times\Bbb P^i\times\{\bold 0\}.
$$

\vskip .1in

\noindent c) there exists an open neighborhood $\Cal W$ of $\bold p$ in $\Cal U$ such that, for all $i$ such that $0\leqslant i\leqslant n$,
$$
\Bbb P(T^*_{S_\alpha}\Cal U)\ \cap \ \Big(\big(\Cal W\cap V(z_0-p_0, \dots, z_{i-1}-p_{i-1})\big)\times\Bbb P^i\times\{\bold 0\}\Big) \ \subseteq \ \{\bold p\}\times\Bbb P^i\times\{\bold 0\}.
$$
(When $i=0$, we mean that $\Bbb P(T^*_{S_\alpha}\Cal U)\ \cap \ \big(\Cal W\times\Bbb P^0\times\{\bold 0\}\big) \ \subseteq \ \{\bold p\}\times\Bbb P^0\times\{\bold 0\}$.)

\vskip .25in

Moreover, we obtain equivalent statements if we replace the condition that $0\leqslant i\leqslant n$ in b) and c) by the condition that $0\leqslant i\leqslant -1+\dm S_\alpha$.
}

\vskip .2in

\noindent{\it Proof}. That a) and b) are equivalent is simply a translation of transversality into conormal terms. Certainly, c) implies b). That we obtain equivalent statements if we replace  $0\leqslant i\leqslant n$ in b) and c) by the condition that $0\leqslant i\leqslant -1+\dm S_\alpha$ is also trivial. The surprising implication is that b) implies the seemingly stronger c).

\vskip .1in

Assume that b) holds. Suppose that we had the germ of a complex analytic curve $\bold q(t):= (\bold x(t), [\omega(t)])\in \Cal U\times \Bbb P^n$ such that $\bold x(0)=\bold p$ and so that, for small $t\neq 0$, $\bold q(t)$ is in 
$$
\Bbb P(T^*_{S_\alpha}\Cal U)\ \cap \ \big( V(z_0-p_0, \dots, z_{k-1}-p_{k-1})\times\Bbb P^k\times\{\bold 0\}\big).
$$
Then, for $t\neq 0$, $\bold x(t) = (\bold 0, x_k(t), \dots, x_n(t))$, $\omega(t) = \omega_0(t)dz_0 +\dots \omega_k(t)dz_k$, $\bold x^\prime(t)\in T_{\bold x(t)}S_\alpha$, and so $\omega_k(t)x^\prime_k(t)\equiv 0$. Thus, either $\omega_k(t)\equiv 0$ or $x_k(t)\equiv p_k$. 

Letting  $Y_k:= V(z_0-p_0, \dots, z_k-p_k)\times\Bbb P^k\times\{\bold 0\}$, the above paragraph shows that, as germs over $\bold p$,
$$
\Bbb P(T^*_{S_\alpha}\Cal U)\ \cap \ \big(V(z_0-p_0, \dots, z_{k-1}-p_{k-1})\times\Bbb P^k\times\{\bold 0\}\big) \ \subseteq \ Y_{k-1}\ \cup \ Y_k,
$$
where $Y_{-1}:=\emptyset$. Now, b) yields a contradiction.\qed

\vskip .3in

\noindent{\bf Corollary 4.4}. {\it Suppose that $\Cal S$ is a Whitney a) partition of $X$. Then the coordinates $\bold z$ are essentially transverse to $\Cal S$ at $\bold p\in X$ if and only if, for all $S_\alpha\in\Cal S$, for all $i$ such that $0\leqslant i\leqslant n$, 
$$\dm_\bold p \left(V(z_0-p_0, \dots, z_{i-1}-p_{i-1})\ \cap\ \nu\big(\Bbb P(\con)\ \cap\ (\Cal U\times \Bbb P^i\times\{\bold 0\})\big)\right) \ \leqslant \ 0.
$$

In particular, if $\bold z$ is essentially transverse to $\Cal S$ at $\bold p$, then $\bold z$ is essentially transverse to $\Cal S$ at all points near $\bold p$.
}

\vskip .2in

\noindent{\it Proof}. This follows from 4.3.c and the fact that the Whitney a) condition implies that $\dsize\bigcup_{\alpha}T^*_{S_\alpha}\Cal U=\bigcup_{\alpha}\con$.\qed

\vskip .3in

As our final result of this section, we wish to characterize Thom's $a_f$ condition in terms of exceptional divisors; this  is essentially Proposition 4.3 from [{\bf M9}]. Recall the definition of the relative conormal space from 3.1.

\vskip .2in

\noindent{\bf Definition 4.5}. Let $M$ and $N$ be analytic submanifolds of $X$ such that
$f$ has constant rank on $N$. Then, the pair $(M, N)$ {\it satisfies Thom's $a_f$ condition at a point $\bold
x\in N$} if and only if we have the containment
$\left(\overline{ T^*_{f_{|_M}}\Cal U}\right)_\bold x \ \subseteq \ 
\Big(T^*_{f_{|_N}}\Cal U\Big)_\bold x$ of fibres over
$\bold x$.

In particular, if $f$ is, in fact, constant on $N$, then the pair $(M, N)$ satisfies Thom's $a_f$ condition at a
point $\bold x\in N$ if and only if we have the containment
$\left(\overline{ T^*_{f_{|_M}}\Cal U}\right)_\bold x \ \subseteq \ 
\Big(T^*_{{}_N}\Cal U\Big)_\bold x$ of fibres over
$\bold x$. If $f$ is constant on both $M$ and $N$, then the $a_f$ condition reduces to Whitney's condition a).

\vskip .3in

\noindent{\it Remark 4.6}. Note that, if $Y$ is an analytic subspace of $X$, and $M$ is an open dense subset of $Y_{\operatorname{reg}}$, then $\overline{T^*_{{}_{Y_{\operatorname{reg}}}}\Cal U}= \overline{T^*_{{}_{M}}\Cal U}$ and, for every submanifold $N\subseteq \Cal U$ and every $\bold x\in N$, $(Y_{\operatorname{reg}}, N)$ satisfies Whitney's condition a)  at $\bold x$ if and only if $(M, N)$ satisfies Whitney's condition a) at $\bold x$; moreover, if $f$ is not constant on any irreducible component of $Y$, then  $(Y_{\operatorname{reg}}, N)$ satisfies Thom's $a_f$ condition  at $\bold x$ if and only if $(M, N)$ satisfies Thom's $a_f$ condition at $\bold x$. Thus, below, it will suffice to work with $Y_{\operatorname{reg}}$ everywhere, instead of the seemingly more general $M$.

The notion of Thom's $a_f$ condition that we use below (and above) is slightly more general than is sometimes the case; we do not require the rank of $f$ to be constant on the bigger stratum. If we were to require the rank of $f$ to be constant on the bigger stratum, then we would be forced to write the more cumbersome `` $(Y_{\operatorname{reg}}-\Sigma(f_{|_{Y_{\operatorname{reg}}}}), N)$ satisfies Thom's $a_f$ condition at $\bold x$''.  Moreover, if $C$ is a component of $Y$ on which $f$ is constant, then saying that $(Y_{\operatorname{reg}}, N)$ satisfies the $a_f$ condition implies that $(C_{\operatorname{reg}}, N)$ satisfies Whitney's condition a); the condition that $(Y_{\operatorname{reg}}-\Sigma(f_{|_{Y_{\operatorname{reg}}}}), N)$ satisfies Thom's $a_f$ condition would ignore what happens on a component such as $C$.

\vskip .3in

\noindent{\bf Definition 4.7}. An {\it $a_f$ partition of $X$ with respect to $\Fdot$} (or, an {\it  $a_{{}_{f, \Fdot}}$ partition of $X$}) is an $\Fdot$-partition, $\Cal S$, of $X$ such that $V(f)$ is a union is strata (the {\it $V(f)$ strata}) and such that for every $\Fdot$-visible $S_\alpha$ and every $V(f)$ stratum $S_\beta$, the pair $(S_\alpha, S_\beta)$ satisfies the $a_f$ condition.

\vskip .3in

The following theorem looks like a significant improvement of Proposition 4.3 of [{\bf M9}]; however, the proof is essentially the same. In Section 6, this theorem will enable us to link Thom's $a_f$ condition with the vanishing cycles along $f$.

\vskip .3in

\noindent{\bf Theorem 4.8}. {\it Suppose that $Y$ is an analytic subset of $X$.   Let $E$ denote the
exceptional divisor in
$\operatorname{Bl}_{\imdf}\overline{T^*_{{}_{Y_{\operatorname{reg}}}}\Cal U}\ \subseteq\ \Cal U\times\Bbb
C^{n+1}\times\Bbb P^n$. Suppose that $N\subseteq X$ is a complex analytic submanifold of $\Cal U$ and that $\bold x\in N$. 

Then,  $(Y_{\operatorname{reg}}, N)$ satisfies Thom's $a_f$ condition at $\bold x$ if and only if 

\vskip .1in

\noindent i)\hskip .2in $(Y_{\operatorname{reg}}, N)$ satisfies Whitney's condition a) at $\bold x$;

\vskip .1in

\noindent ii)\hskip .19in $d_\bold x\tilde f\in\big(T^*_{{}_{N}}\Cal U\big)_\bold x$;

\vskip .1in

\noindent iii)\hskip .17in there is the containment of fibres
above
$\bold x$ given by $\big(\pi(E)\big)_{\bold x}\subseteq \big(\Bbb P(T^*_{{}_N}\Cal U)\big)_\bold x$.
}

\vskip .2in

\noindent{\it Proof}. We first make a few simple observations.

\vskip .2in

\noindent $\bullet$\hskip .2in Suppose $(Y_{\operatorname{reg}}, N)$ satisfies Thom's $a_f$ condition at $\bold x$. 

\vskip .1in

Then certainly i) and ii) follow; for we can use as relative conormal covectors every conormal to $Y_{\operatorname{reg}}$ and every covector of the form $d_\bold y\tilde f$, where $\bold y\in Y_{\operatorname{reg}}$. In addition, we proved in Proposition 4.3 of [{\bf M9}] that the $a_f$ condition implies that $$\big(\pi(E)\big)_{\bold x}\subseteq \big(\Bbb P(T^*_{{}_N}\Cal U)\big)_\bold x.$$

\vskip .1in

This proves one direction of the theorem.

\vskip .2in

\noindent $\bullet$\hskip .2in Now suppose that i), ii), and iii) hold.

\vskip .1in

If  $C$ is a component of $Y$ on which $f$ is constant, then limiting relative conormals to $C_{\operatorname{reg}}$ are the same as limiting (absolute) conormals to $C_{\operatorname{reg}}$. As we are assuming that Whitney's condition a) holds, it follows that $(C_{\operatorname{reg}}, N)$ satisfies that $a_f$ at $\bold x$.

For components $C$ of $Y$ on which $f$ is non-constant, the proof of  Proposition 4.3 of [{\bf M9}] shows that  $(C_{\operatorname{reg}}, N)$ satisfies that $a_f$ at $\bold x$.\qed

\vskip .3in

\noindent\S5. {\bf Characteristic Polar Complexes, Modules, and Isolating Coordinates}

\vskip .1in

In this section, we will define one our primary objects of study -- {\it the characteristic polar modules} -- and relate them to polar varieties, graded, enriched characteristic cycles, and iterated vanishing and nearby cycles. This will require us to investigate how generic a linear choice of coordinates must be in order to produce nice results.
It is important throughout this section, and throughout the remainder of this paper, that our notion of {\it isolating coordinates} is an {\bf effective} notion of ``generic''; that is, in many situations, one can determine fairly easily whether the coordinates are isolating. Moreover, if the coordinates are isolating at a given point, then they are isolating at all nearby points -- there is no need to re-choose  the coordinates at each point. These properties are important, since our goal is to produce effectively calculable data that one can associate to a singularity.

\vskip .1in

We continue with all of our previous notation, and also introduce new notation. 

\vskip .1in

We let $Y$ be a new complex analytic subspace of $\Cal U$ and we let $\Adot$ be a bounded, constructible complex of $R$-modules on $Y$. Let  $\Cal R:=\{R_\beta\}$ be an $\Adot$-partition of $Y$. Our reason for introducing these new objects is that, later, we will return to the setting of the previous sections by considering the special case where $Y=V(f)=X\cap V(\tilde f)$ and $\Adot=\phi_f[-1]\Fdot$.

\vskip .1in

Fix a point $\bold p=(p_0, \dots, p_n)\in Y$, and let $d:=\dim_\bold p (\supp\Adot)$. 

\vskip .1in

Recall that we use $\bold z=(z_0, \dots,
z_n)$ to denote coordinates on $\Cal U$. Below, when the context makes the domains clear, we shall not distinguish in the
notation between the coordinate functions $z_i$ and their restrictions to various subspaces.  As we will be
projectivizing conormal varieties, we  assume that $\operatorname{codim}_{{}_{\Cal U}}Y\geqslant 1$.

Throughout the remainder of this paper, it will be convenient to adopt the standard convention that, when $m=0$,  $V(z_0-p_0, \dots, z_{m-1}-p_{m-1}) = \Cal U$.

\vskip .3in

This section contains a number of technical results. However, the reader should note that the main points of this section are:

\vskip .1in

\noindent $\bullet$\hskip .2in In Definition 5.1, we define the {\it $j$-th characteristic polar module in degree $k$ of $\Adot$ with respect to the coordinates $\bold z$ at the point $\bold p$}  to be 
$${}^k\gamma^j_{{}_{\Adot, \bold z}}(\bold p)\ :=\ H^k\big( \phi_{z_j-p_j}[-1]\psi_{z_{j-1}-p_{j-1}}[-1]\dots\psi_{z_0-p_0}[-1]\Adot\big)_\bold p.$$

\vskip .1in

\noindent$\bullet$\hskip .2in The coordinates $\bold z$ are {\it $\Adot$-isolating at $\bold p$} provided that the support of $$ \phi_{z_j-p_j}[-1]\psi_{z_{j-1}-p_{j-1}}[-1]\dots\psi_{z_0-p_0}[-1]\Adot$$ is at most $0$-dimensional at $\bold p$ for all $j$. Despite the fact that this condition looks fairly unmanageable, in Theorem 5.10, we prove that this condition has a nice interpretation in terms of intersections with $\gecc^\bullet(\Adot)$. Moreover, under the assumption that the coordinates are $\Adot$-isolating, general results on nearby cycles, vanishing cycles, and perverse cohomology allow us to conclude Theorem 5.18: a result which yields chain complexes, containing the ${}^k\gamma^j_{{}_{\Adot, \bold z}}(\bold x)$, whose cohomology is isomorphic to the stalk cohomology of the perverse cohomology of $\Adot$ in each degree. This result is important because the modules appearing in these chain complexes can be calculated; see the following paragraph.

\vskip .1in

\noindent$\bullet$\hskip .2in In Definition 5.20, under the assumption that the coordinates are $\Adot$-isolating, we define the {\it enriched $j$-th characteristic polar cycle in degree $k$ of $\Adot$ with respect to the coordinates $\bold z$} to be 
$$
{}^k\Gamma^j_{{}_{\Adot, \bold z}}:= \nu_*\Big(\Bbb P(\gecc^k(\Adot)) \ \odot\ \Cal W\times\Bbb
P^j\times\{\bold 0\}\Big)
$$
and we show in Theorem 5.23 that ${}^k\gamma^j_{{}_{\Adot, \bold z}}(\bold x)$ can, in fact, be calculated by 
$$ {}^k\gamma^j_{{}_{\Adot, \bold z}}(\bold p) \ \cong \ \big({}^k\Gamma^j_{{}_{\Adot, \bold z}} \ \odot \ V(z_0-p_0, \dots, z_{j-1}-p_{j-1})\big)_\bold p.$$

\vskip .3in

\noindent{\bf Definition 5.1}. For all $j$ such that $0\leqslant j\leqslant n$,  for all $\bold a\in\Bbb C^{j+1}$, we define the {\it $j$-th characteristic polar complex of $\Adot$ with respect to $\bold z$ at $\bold a$} to be
$$\Phi^j_{{}_{\Adot, \bold z}}(\bold a) \ := \ \phi_{z_j-a_j}[-1]\psi_{z_{j-1}-a_{j-1}}[-1]\dots\psi_{z_0-a_0}[-1]\Adot.$$

When $j=0$, we mean that $\Phi^0_{{}_{\Adot, \bold z}}(a_0) \ := \ \phi_{z_0-a_0}[-1]\Adot$. For $\bold x\in\Bbb C^{n+1}$, it is convenient to define $\Phi^j_{{}_{\Adot, \bold z}}(\bold x) \ := \ \Phi^j_{{}_{\Adot, \bold z}}(x_0, x_1, \dots, x_j)$.

\vskip .1in

For all $\bold x\in Y$, we define the {\it $j$-th characteristic polar module in degree $k$ of $\Adot$ with respect to $\bold z$ at $\bold x$}, ${}^k\gamma^j_{{}_{\Adot, \bold z}}(\bold x)$,  to be the degree $k$ stalk cohomology at $\bold x$ of $\Phi^j_{{}_{\Adot, \bold z}}(\bold x)$, i.e., 
$${}^k\gamma^j_{{}_{\Adot, \bold z}}(\bold x)\ :=\ H^k\big(\Phi^j_{{}_{\Adot, \bold z}}(\bold x)\big)_\bold x.$$

The {\it support of ${}^k\gamma^j_{{}_{\Adot, \bold z}}$} is, naturally, defined to be the closure of the set $\{\bold x\in Y\ |\ {}^k\gamma^j_{{}_{\Adot, \bold z}}(\bold x)\neq 0\}$. We denote this support by $\supp({}^k\gamma^j_{{}_{\Adot, \bold z}})$.

We define the {\it support of ${}^\bullet\gamma^j_{{}_{\Adot, \bold z}}$} to be the closure of the set $$\{\bold x\in Y\ |\text{ there exists }k\text{ such that }\ {}^k\gamma^j_{{}_{\Adot, \bold z}}(\bold x)\neq 0\}.$$
We denote this support by $\supp({}^\bullet\gamma^j_{{}_{\Adot, \bold z}})$.
By boundedness, there are a finite number of $k$ such that ${}^k\gamma^j_{{}_{\Adot, \bold z}}\neq 0$; hence,
$$
\supp({}^\bullet\gamma^j_{{}_{\Adot, \bold z}})\ =\ \bigcup_k \supp({}^k\gamma^j_{{}_{\Adot, \bold z}}).
$$
\vskip .1in

If the base ring $R$ is a domain, we define the {\it  $j$-th polar Euler number of $\Adot$ with respect to $\bold z$ at $\bold x$}
to be 
$$ {}^{\operatorname{ord}}\gamma^j_{{}_{\Adot, \bold z}}(\bold x) = \sum_k(-1)^k\operatorname{rk}\big({}^{k}\gamma^j_{{}_{\Adot,
\bold z}}(\bold x)\big).
$$

\vskip .2in

We make the following trivial observation.

\vskip .2in

\noindent{\bf Proposition 5.2}. {\it Let $\tilde{\bold z}$ denote the ``rotated'' coordinate system $(z_1, z_2, \dots, z_n, z_0)$. Fix $x_0$. 

For all $j$ such that $1\leqslant j\leqslant n$, for all $\bold x\in V(z_0-x_0)$,
$$\Phi^{j-1}_{\psi_{z_0-x_0}[-1]\Adot, \tilde{\bold z}}(\bold x)=\Phi^j_{{}_{\Adot, \bold z}}(\bold x),$$
and, for all $k$,
$$
{}^k\gamma^{j-1}_{{}_{\psi_{z_0-x_0}[-1]\Adot, \tilde{\bold z}}}(\bold x) \ = \ {}^k\gamma^j_{{}_{\Adot, \bold z}}(\bold x).
$$
}

\vskip .2in

\noindent{\it Proof}. This is immediate from the definitions.\qed

\vskip .3in

Only slightly less trivial is the following.

\vskip .2in

\noindent{\bf Proposition 5.3}. {\it 
If $R$ is a domain, then, for all $m$ such that $0\leqslant m\leqslant n$, the Euler characteristic of the stalk of $\Adot$ at $\bold x$ is given by
$$\chi(\Adot)_\bold x = (-1)^{m+1}\chi(\psi_{z_m-x_m}[-1]\dots\psi_{z_0-x_0}[-1]\Adot)_\bold x \ + \ \sum_{j\leqslant m}(-1)^j\big({}^{\operatorname{ord}}\gamma^j_{{}_{\Adot, \bold
z}}(\bold x)\big).$$
}

\vskip .2in

\noindent{\it Proof}. This follows immediately by inductively applying the fact that Euler characteristics are
additive over long exact sequences to the long exact sequences on stalk cohomology, at $\bold x$, which come from the
distinguished triangles
$$
\Adot_{|_{V(z_j-x_j)}}[-1] \ \rightarrow \ \psi_{z_j-x_j}[-1]\Adot \ \rightarrow \ \phi_{z_j-x_j}[-1]\Adot\
@>[1]>>\Adot_{|_{V(z_j-x_j)}}[-1].\qed
$$

\vskip .3in

While it is clear that $$\supp(\Phi^j_{{}_{\Adot, \bold z}}(\bold x))\subseteq V(z_0-x_0, \dots, z_j-x_j)\cap \supp({}^\bullet\gamma^j_{{}_{\Adot, \bold z}}),$$ the reverse inclusion need not hold. In order to determine a more precise relationship between $\supp(\Phi^j_{{}_{\Adot, \bold z}}(\bold x))$ and \hbox{$V(z_0-x_0, \dots, z_j-x_j)\cap \supp({}^\bullet\gamma^j_{{}_{\Adot, \bold z}})$,} we need a lemma and a proposition.

\vskip .3in

\noindent{\bf Lemma 5.4}. {\it Let $h:Y\rightarrow\Bbb C$ be an analytic function, and let $\tilde h:\Cal U\rightarrow\Bbb C$ be a local extension of $h$ near $\bold p$. Suppose that $\bold p\not\in\supp\phi_{h-h(\bold p)}[-1]\Adot$. Then,  for all $k$, there
is an equality of fibres over
$\bold p$ given by
$$
\big|\big(T^*_{{}_{h, \Adot}}\Cal U\big)^k\big|_\bold p \ = \ |\gecc^k(\Adot)|_\bold p \ + \ <d_\bold p\tilde h>,
$$ where $<d_\bold p\tilde h>$ denotes the linear subspace of all scalar multiples of $d_\bold p\tilde h$. }

\vskip .2in

\noindent{\it Proof}.  By the definition of $T^*_{h_{|_{R_\beta}}}\Cal U$, if $\bold p\in R_\beta$ and $d_\bold p\tilde h\not\in T^*_{{}_{R_\beta}}\Cal U$, then $$\big(T^*_{h_{|_{R_\beta}}}\Cal U\big)_\bold p = \big(T^*_{{}_{R_\beta}}\Cal U\big)_\bold p+\big<d_\bold p\tilde h\big>.$$
We need to deal with possible strata $R_\beta$ such that $\bold p\in\overline{R_\beta}-R_\beta$.

\vskip .1in

 Fix $k$. Since $\bold p\not\in\supp\phi_{h-h(\bold p)}[-1]\Adot$, Theorem 3.4 implies that $d_\bold p\tilde h\not\in |\gecc^k(\Adot)|_\bold p$. Suppose that $\bold p\in\overline{R_\beta}$ and $\overline{T^*_{R_\beta}\Cal U}$ is a component of $|\gecc^k(\Adot)|$; then, $d_\bold p\tilde h\not\in \big(\overline{T^*_{{}_{R_\beta}}\Cal U}\big)_\bold p$.

By 3.4, the set $\eta\big(|\gecc^\bullet(\Adot)|\cap\operatorname{im}d\tilde h\big)$ is closed, and so,  for all $\bold x$ near $\bold p$, $d_\bold x\tilde h\not\in |\gecc^k(\Adot)|_\bold x$ . Therefore, if , we may take $\bold x_i\in R_\beta$ such that $\bold x_i\rightarrow \bold p$ and
$$\big(T^*_{h_{|_{R_\beta}}}\Cal U\big)_{\bold x_i} = \big(T^*_{{}_{R_\beta}}\Cal U\big)_{\bold x_i}+\big<d_{\bold x_i}\tilde h\big>.$$

As $d_\bold p\tilde h\not\in \big(\overline{T^*_{{}_{R_\beta}}\Cal U}\big)_\bold p$, we conclude that the limits behave ``nicely'', and thus $$\big(\overline{T^*_{h_{|_{R_\beta}}}\Cal U}\big)_\bold p = \big(\overline{T^*_{{}_{R_\beta}}\Cal U}\big)_\bold p+\big<d_\bold p\tilde h\big>.\qed$$

\vskip .3in

\noindent{\bf Definition 5.5}  For all $m$ such that  $0\leqslant m\leqslant n$, we
define {\it $\Theta^m_{{}_{\Adot, \bold z}}$ } by 
$$
 \Theta^m_{{}_{\Adot, \bold z}} \ := \  \nu\Big(|\Bbb P(\gecc^\bullet(\Adot))| \ \cap \ (\Cal U\times\Bbb
P^m\times\{\bold 0\})\Big),
$$
 and we define {\it $\Gamma^m_{{}_{\Adot, \bold z}}$} to be the union of the $m$-dimensional components of
$\Theta^m_{{}_{\Adot, \bold z}}$. We refer to $\Gamma^m_{{}_{\Adot, \bold z}}$ as the {\it  $m$-dimensional characteristic polar variety of
$\Adot$ with respect to $\bold z$}.

\vskip .3in

\noindent{\it Remark 5.6}. Note that each $ \Theta^m_{{}_{\Adot, \bold z}}$ is closed. In addition, if  $R_\beta$ is a maximal stratum of $\supp\Adot$, and $d_\beta:=\dim R_\beta$, then the fibres of $|\Bbb P(\gecc^\bullet(\Adot))|$ over $R_\beta$ will be of dimension $n-d_\beta$, and so these fibres must intersect $\Bbb P^{d_\beta}\times\{\bold 0\}$; therefore, $\overline{R_\beta}\subseteq\Theta^{d_\beta}_{{}_{\Adot, \bold z}}$.

It follows that, if $d:=\dim_\bold p\supp(\Adot)$, then, for all $m\geqslant d$, as germs of sets at $\bold p$, $\Theta^m_{{}_{\Adot, \bold z}} = \supp(\Adot)$.

\vskip .1in

The sets $\Theta^m_{{}_{\Adot, \bold z}}$ are very closely related to the absolute polar varieties of L\^e and Teissier. The set $\nu\big(\overline{T^*_{R_\beta}\Cal U} \ \cap \ (\Cal U\times\Bbb P^m\times\{\bold 0\})\big)$ consists of the closure of the critical locus of the map $(z_0, \dots, z_m)_{|_{R_\beta}}$ together with some possible ``degenerate'' critical points on smaller strata. By Whitney's condition a), these degenerate points will be critical points of $(z_0, \dots, z_m)$ restricted to smaller strata. Hence, we have the containments
$$
\bigcup\Sb\Adot-\text{visible } \\ R_\beta\endSb\overline{\operatorname{crit}(z_0, \dots, z_m)_{|_{R_\beta}}} \ \subseteq \ \Theta^m_{{}_{\Adot, \bold z}} \  \subseteq \  \bigcup\Sb\Adot-\text{visible } \\ R_\beta\endSb \bigcup_{R_\gamma\subseteq \overline{R_\beta}}\overline{\operatorname{crit}(z_0, \dots, z_m)_{|_{R_\gamma}}}.
$$
If the coordinates $\bold z$ are sufficiently generic at $\bold p$, and $\dim R_\gamma\geqslant m$, then $\overline{\operatorname{crit}(z_0, \dots, z_m)_{|_{R_\gamma}}}$ is precisely the $m$-dimensional absolute polar variety of $\overline{R_\gamma}$ at $\bold p$.

\vskip .3in

\noindent{\bf Proposition 5.7}. {\it For all $m$ such that $0\leqslant m\leqslant n$, there is an equality of sets given by
$$
\Theta^m_{{}_{\Adot, \bold z}} = \bigcup_{0\leqslant j\leqslant m}\Big(\bigcup_{\bold a\in\Bbb C^{j+1}}\supp\big(\Phi^j_{{}_{\Adot, \bold z}}(\bold a)\big)\Big).
$$
}

\vskip .1in

\noindent{\it Proof}. We proceed by induction on $m$.

\vskip .1in

When $m=0$, the claim is equivalent to
$$
\eta\Big(|\gecc^\bullet(\Adot)| \ \cap \ \operatorname{im}dz_0\Big) = \bigcup_{a_0\in \Bbb C}\supp(\phi_{z_0-a_0}[-1]\Adot),
$$
which follows from the first equality of Theorem 3.4.

\vskip .1in

Now suppose the claim is true for $m-1$, for all complexes and all choices of coordinates, where $1\leqslant m\leqslant n$. We wish to prove the claim for $m$. 

Fix $a_0\in \Bbb C$. It suffices to prove that 
$$
V(z_0-a_0)\cap \Theta^m_{{}_{\Adot, \bold z}} = \supp(\phi_{z_0-a_0}[-1]\Adot) \ \cup \ \bigcup_{1\leqslant j\leqslant m}\Big(\bigcup_{(a_1, \dots, a_j)\in\Bbb C^{j}}\supp\big(\Phi^j_{{}_{\Adot, \bold z}}(a_0, a_1, \dots, a_j)\big)\Big).\tag{$\dagger$}
$$

\vskip .1in

Fix $\bold p\in V(z_0-a_0)\cap Y$. We wish to show that $\bold p$ is in the set on the left side of $(\dagger)$ if and only if $\bold p$ is in the set on the right side of $(\dagger)$. There are two cases.

\vskip .3in

\noindent{\bf case 1}: $\bold p\in \supp(\phi_{z_0-a_0}[-1]\Adot)$.

\vskip .1in

This case follows trivially from the $m=0$ discussion above.

\vskip .3in

\noindent{\bf case 2}: $\bold p\not\in\supp(\phi_{z_0-a_0}[-1]\Adot)$.

\vskip .1in

Recalling Proposition 5.2, we apply the inductive hypothesis to $\psi_{z_0-a_0}[-1]\Adot$ and the rotated coordinates $\tilde{\bold z}$ to obtain that
$$
\nu\Big(|\Bbb P(\gecc^\bullet(\psi_{z_0-a_0}[-1]\Adot))| \ \cap \ (\Cal U\times\{0\}\times\Bbb
P^{m-1}\times\{\bold 0\})\Big) \ =$$
 $$\bigcup_{1\leqslant j\leqslant m}\Big(\bigcup_{(a_1, \dots, a_j)\in\Bbb C^{j}}\supp\big(\Phi^j_{{}_{\Adot, \bold z}}(a_0, a_1, \dots, a_j)\big)\Big).
$$

By Theorem 3.2, $\gecc^\bullet(\psi_{z_0-a_0}[-1]\Adot)  \ = \ \big(T^*_{{}_{z_0, \Adot}}\Cal U\big)^\bullet\odot(V(z_0-a_0)\times\Bbb C^{n+1})$. As we are in the case where $\bold p\not\in\supp\phi_{z_0-a_0}[-1]\Adot$, Lemma 5.4 tells us that $$\big|\big(T^*_{{}_{z_0, \Adot}}\Cal U\big)^\bullet\big|_\bold p \ = \ |\gecc^\bullet(\Adot)|_\bold p \ + \ <d_\bold p  z_0>.$$

It is now trivial to show that: $\bold p\in V(z_0-a_0)\cap \nu\Big(|\Bbb P(\gecc^\bullet(\Adot))| \ \cap \ (\Cal U\times\Bbb
P^m\times\{\bold 0\})\Big)$ if and only if 
$$
\bold p\in \nu\Big(|\Bbb P(\gecc^\bullet(\psi_{z_0-a_0}[-1]\Adot))| \ \cap \ (\Cal U\times\{0\}\times\Bbb
P^{m-1}\times\{\bold 0\})\Big),
$$
if and only if
$$
\bold p\in \bigcup_{1\leqslant j\leqslant m}\Big(\bigcup_{(a_1, \dots, a_j)\in\Bbb C^{j}}\supp\big(\Phi^j_{{}_{\Adot, \bold z}}(a_0, a_1, \dots, a_j)\big)\Big),$$
i.e., since we are in case 2, that $\bold p$ is in the set on the left side of $(\dagger)$ if and only if $\bold p$ is in the set on the right side of $(\dagger)$.\qed

\vskip .3in

\noindent{\bf Corollary 5.8}. {\it For all $m$ such that $0\leqslant m\leqslant n$, there is an equality of sets given by
$$
\Theta^m_{{}_{\Adot, \bold z}} = \bigcup_{0\leqslant j\leqslant m}\supp({}^\bullet\gamma^j_{{}_{\Adot, \bold z}}).
$$
}

\vskip .2in

\noindent{\it Proof}. This is actually a point-set topology proof. Let $E^j_{{}_{\Adot, \bold z}}(\bold a)$ denote the set of those $\bold x\in Y$ such that the stalk of $\Phi^j_{{}_{\Adot, \bold z}}(\bold a)$ at $\bold x$ is not zero. Hence, $\supp(\Phi^j_{{}_{\Adot, \bold z}}(\bold a))=\overline{E^j_{{}_{\Adot, \bold z}}(\bold a)}$, and $\supp({}^\bullet\gamma^j_{{}_{\Adot, \bold z}}) = \overline{\bigcup_{\bold a\in\Bbb C^{j+1}}E^j_{{}_{\Adot, \bold z}}(\bold a)}$.

\vskip .1in

Now, Proposition 5.7 tells us that $\dsize\bigcup_{0\leqslant j\leqslant m}\Big(\bigcup_{\bold a\in\Bbb C^{j+1}}\overline{E^j_{{}_{\Adot, \bold z}}(\bold a)}\Big)$
is closed. It is a trivial topology proof to show then that 
$$
\bigcup_{0\leqslant j\leqslant m}\Big(\bigcup_{\bold a\in\Bbb C^{j+1}}\overline{E^j_{{}_{\Adot, \bold z}}(\bold a)}\Big) \ = \ \bigcup_{0\leqslant j\leqslant m}\overline{\Big(\bigcup_{\bold a\in\Bbb C^{j+1}}E^j_{{}_{\Adot, \bold z}}(\bold a)\Big)}.\qed
$$

\vskip .3in

\noindent{\bf Lemma 5.9}. {\it Fix $m$. If, for all $j\leqslant m$, $\dim(\Cal W\cap \Theta^j_{{}_{\Adot, \bold z}})\leqslant j$, then, for all $j\leqslant m$, $|\Bbb P(\gecc^\bullet(\Adot))|$ properly intersects $\Cal W\times \Bbb P^j\times\{\bold 0\}$ inside $\Cal U\times\Bbb P^n$.

Moreover, whenever $|\Bbb P(\gecc^\bullet(\Adot))|$ properly intersects $\Cal W\times \Bbb P^j\times\{\bold 0\}$ inside $\Cal U\times\Bbb P^n$ for all $j\leqslant m$, there is an equality of sets, inside of $\Cal W$, given by $\Theta^m_{{}_{\Adot, \bold z}}=\bigcup_{j\leqslant m} \Gamma^j_{{}_{\Adot, \bold z}}$.

Finally, fix $j$, where $0\leqslant j\leqslant n$. Suppose that the analytic set
$|\Bbb P(\gecc^\bullet(\Adot))|$ properly intersects $\Cal W\times\Bbb P^j\times\{\bold 0\}$ inside $\Cal U\times \Bbb
P^n$. Then,  the graded enriched cycle $\Bbb P(\gecc^\bullet(\Adot))$ properly intersects $\Cal W\times\Bbb P^j\times\{\bold 0\}$ inside $\Cal U\times \Bbb P^n$, and, inside of $\Cal W$,
$$
\Gamma^j_{{}_{\Adot, \bold z}} \ = \ \Big|\nu_*\Big(\Bbb P(\gecc^\bullet(\Adot)) \ \odot\ \Cal W\times\Bbb
P^j\times\{\bold 0\}\Big)\Big|,
$$
where $\nu_*$ denotes the proper push-forward.
}

\vskip .2in

\noindent{\it Proof}. Recall that we use $[w_0:w_1:\dots:w_n]$ for homogeneous coordinates on $\Bbb P^n$. 

\vskip .1in

Assume that, for all $j\leqslant m$, $\dim(\Cal W \cap \Theta^j_{{}_{\Adot, \bold z}})\leqslant j$, and suppose that $E$ is a component of $|\Bbb P(\gecc^\bullet(\Adot))|\cap (\Cal W\times \Bbb P^j\times\{\bold 0\})|$ such that $\dim E\geqslant j+1$, i.e., $E$ is non-proper component of the intersection. As $\eta(E)\subseteq \Theta^j_{{}_{\Adot, \bold z}}$, $\dim(\eta(E))\leqslant j$, by our hypothesis. Therefore, the generic fibre of $E$ over $\eta(E)$ has dimension at least $1$. Consequently, the generic fibre of $E$ intersects the copy of $\Bbb P^{n-1}$ given by  $V(w_j)$. Hence, $E\cap V(w_j)\subseteq |\Bbb P(\gecc^\bullet(\Adot))|\cap (\Cal W\times \Bbb P^{j-1}\times\{\bold 0\})$ and $\dim(E\cap V(w_j))\geqslant j$.

Thus, by induction, we arrive at the fact that $|\Bbb P(\gecc^\bullet(\Adot))|\cap (\Cal W\times \Bbb P^0\times\{\bold 0\})$ has a component with a generic fibre of dimension at least $1$. This contradiction proves that, for all $j\leqslant m$, $|\Bbb P(\gecc^\bullet(\Adot))|$ properly intersects $\Cal W\times \Bbb P^j\times\{\bold 0\}$ inside $\Cal U\times\Bbb P^n$.

\vskip .2in

Now, assume that  $|\Bbb P(\gecc^\bullet(\Adot))|$ properly intersects $\Cal W\times \Bbb P^j\times\{\bold 0\}$ inside $\Cal U\times\Bbb P^n$ for all $j\leqslant m$. 

If $j\leqslant m$, then $\Theta^j_{{}_{\Adot, \bold z}}\subseteq \Theta^m_{{}_{\Adot, \bold z}}$. As $\Gamma^j_{{}_{\Adot, \bold z}}\subseteq \Theta^j_{{}_{\Adot, \bold z}}$, it follows that $\bigcup_{j\leqslant m} \Gamma^j_{{}_{\Adot, \bold z}}\subseteq \Theta^m_{{}_{\Adot, \bold z}}$. We need to show the reverse containment.

Suppose that $C$ is an irreducible component of $ \Theta^m_{{}_{\Adot, \bold z}}$ and that $E$ is an irreducible component of $|\Bbb P(\gecc^\bullet(\Adot))|\cap (\Cal W\times \Bbb P^m\times\{\bold 0\})$ such that $\eta(E)=C$.  By hypothesis, $\dim E=m$, and so $\dim C\leqslant m$. 

Suppose that $\dim C\leqslant m-1$, i.e., that $C$ is a component of $ \Theta^m_{{}_{\Adot, \bold z}}$, but not a component of $ \Gamma^m_{{}_{\Adot, \bold z}}$. Then, the generic fibre of $E$ must be at least $1$-dimensional. Hence, $V(w_m)$ intersects the generic fibres of $E$, and so $C=\eta(E\cap V(w_m))\subseteq \Theta^{m-1}_{{}_{\Adot, \bold z}}$. Proceeding inductively, we conclude that: if $C$ is a $j$-dimensional irreducible component of $ \Theta^m_{{}_{\Adot, \bold z}}$, then $C$ is a component of $ \Theta^j_{{}_{\Adot, \bold z}}$ and, therefore, a component of $ \Gamma^j_{{}_{\Adot, \bold z}}$. This proves that $\Theta^m_{{}_{\Adot, \bold z}}=\bigcup_{j\leqslant m} \Gamma^j_{{}_{\Adot, \bold z}}$.

\vskip .2in

Finally, suppose, for a fixed $j$, that the analytic set
$|\Bbb P(\gecc^\bullet(\Adot))|$ properly intersects $\Cal W\times\Bbb P^j\times\{\bold 0\}$ inside $\Cal U\times \Bbb
P^n$.

As every component of $\Bbb P(\gecc^\bullet(\Adot))$ is purely $n$-dimensional, the components of $\Bbb P(\gecc^\bullet(\Adot))$ are the same as the components of the underlying set $|\Bbb P(\gecc^\bullet(\Adot))|$. Thus, $\Bbb P(\gecc^\bullet(\Adot))$ properly intersects $\Cal W\times\Bbb P^j\times\{\bold 0\}$ inside $\Cal U\times \Bbb P^n$, and so the dimension of every component of $\Bbb P(\gecc^\bullet(\Adot)) \  \odot\  \Cal W\times\Bbb
P^j\times\{\bold 0\}$ is equal to $j$. Therefore, by definition of the proper push-forward, the components of $\Big|\nu_*\Big(\Bbb P(\gecc^\bullet(\Adot)) \ \odot\ \Cal W\times\Bbb
P^j\times\{\bold 0\}\Big)\Big|$ are precisely the $j$-dimensional components of $\nu\Big(|\Bbb P(\gecc^\bullet(\Adot))| \ \cap \ (\Cal W\times\Bbb
P^j\times\{\bold 0\})\Big)$.\qed

\vskip .3in

\noindent{\bf Theorem 5.10}. {\it Let $\bold p\in Y$, and fix $m$ such that $0\leqslant m\leqslant n$. Then, the following are equivalent:

\vskip .1in

\noindent a) for all $j$ such that $0\leqslant j\leqslant m$, $\dim_\bold p\supp(\Phi^j_{{}_{\Adot, \bold z}}(\bold p))\leqslant 0$;

\vskip .1in

\noindent b) for all $j$ such that $0\leqslant j\leqslant m$, $\dim_\bold p \Big(V(z_0-p_0, \dots, z_{j-1}-p_{j-1})\cap\supp({}^\bullet\gamma^j_{{}_{\Adot, \bold z}})\Big)\leqslant 0$;

\vskip .1in

\noindent c)  for all $j$ such that $0\leqslant j\leqslant m$,  $\dim_\bold p \Big(V(z_0-p_0, \dots, z_{j-1}-p_{j-1})\cap \Theta^j_{{}_{\Adot, \bold z}}\Big)\leqslant 0$.

\vskip .1in

\noindent d)  for all $j$ such that $0\leqslant j\leqslant m$, there exists an open neighborhood $\Cal W$ of $\bold p$ in $\Cal U$ such that 

$\left|\Bbb P\big(\gecc^\bullet(\Adot)\big)\right|$ properly intersects $\Cal W\times\Bbb P^j\times\{\bold 0\}$ inside $\Cal W\times\Bbb P^n$ and $$\dim_\bold p \Big(V(z_0-p_0, \dots, z_{j-1}-p_{j-1})\cap \left|\nu_*\big(\Bbb P\big(\gecc^\bullet(\Adot)\big)\ \odot\ \Cal W\times\Bbb P^j\times\{\bold 0\}\big)\right|\Big)\leqslant 0.$$

}

\vskip .2in

\noindent{\it Proof}. Given the results of 5.7, 5.8, and 5.9, this is an easy exercise; we leave it to the reader.\qed

\vskip .3in

\noindent{\bf Definition 5.11}.  The coordinates $\bold z=(z_0, \dots, z_n)$ are {\it $\Adot$-isolating at $\bold p$} if and only
if the equivalent conditions of Theorem 5.10 hold for $m=d-1$.

 Here, when $d\leqslant 0$, we mean that there is no condition on the coordinates.

\vskip .1in

Note that this is a condition on only the first $d$ coordinates, and occasionally we will simply say that $(z_0, \dots,
z_{d-1})$ are
$\Adot$-isolating at $\bold p$. Note, also, that if $\bold z$ is $\Adot$-isolating at $\bold p$, then 5.10.c implies that, for all $j$ such that $0\leqslant j\leqslant d-1$, $\dim_\bold p\Theta^j_{{}_{\Adot, \bold z}}\leqslant j$.

\vskip .1in

If $Z\subseteq Y$, then we say that the coordinates $\bold z$ are {\it $\Adot$-isolating on $Z$} if and only if, for all $\bold
x\in Z$, $\bold z$ is $\Adot$-isolating at $\bold x$.

 \vskip .3in

\noindent{\it Remark 5.12}. We could use 5.10.a to define a seemingly more general notion, that of an {\it $\Adot$-isolating sequence of arbitrary functions $(g_0, \dots, g_n)$}, or perhaps a better term would be an {\it $\Adot$-regular sequence}.

However, this notion can be easily recovered from our current set-up by considering the graph map $G:\Cal U\rightarrow\Bbb C^{n+1}\times\Cal U$, given by $G(\bold x) = (g_0(\bold x), \dots, g_n(\bold x), \bold x)$. Use $(u_0,\dots, u_n, z_0, \dots, z_n)$ as coordinates on $\Bbb C^{n+1}\times\Cal U$. Then, the sequence $(g_0, \dots, g_n)$ is $\Adot$-isolating if and only if $(u_0,\dots, u_n)$ is $RG_*\Adot$-isolating.

\vskip .3in

Some immediate properties of $\Adot$-isolating coordinate are:

\vskip .3in

\noindent{\bf Proposition 5.13}. {\it  If $\Cal R$ is Whitney a) and $\bold z$ is essentially transverse to $\Cal R$ at $\bold p$, then $\bold z$ is $\Adot$-isolating at $\bold p$. In particular, $\Adot$-isolating coordinates are generic.

\vskip .1in

If $\bold z$ is $\Adot$-isolating at $\bold p$, there exists an open neighborhood of $\bold p$ on which $\bold z$ is $\Adot$-isolating.

\vskip .1in

The coordinates $\bold z$ are $\Adot$-isolating if and only if $\bold z$ are ${}^\mu\hskip -.02in H^k(\Adot)$-isolating for all $k$, and, in this case, for all $k$ and for all $j$, $0\leqslant j\leqslant n$,
$$
{}^k\gamma^j_{{}_{\Adot, \bold z}}(\bold x)={}^0\gamma^j_{{}_{{}^\mu\hskip -.02in H^k(\Adot), \bold z}}(\bold x).
$$
}

\vskip .1in

\noindent{\it Proof}. The first claim follows from Corollary 4.4.  In the second claim, the existence of the open neighborhood on which $\bold z$ is $\Adot$-isolating follows from the fact that that fibre dimension of $(z_0, \dots, z_{j-1})_{|_{\Theta^j_{{}_{\Adot, \bold z}}}}$ is upper-semicontinuous.

The remainder of the proposition is immediate from the characterization of isolating given in 5.10.a together with the well-known facts that ${}^\mu\hskip -.02in H^k$ commutes with the shifted nearby and vanishing cycles, and that
$$
\supp\Adot = \bigcup_k\supp{}^\mu\hskip -.02in H^k(\Adot).\qed
$$

\vskip .3in

\noindent{\bf Proposition 5.14}. {\it Suppose that the coordinates $\bold z$ are  $\Adot$-isolating at $\bold p$.

\vskip .1in

Then, there exists an open neighborhood $\Cal W$ of $\bold p$ such that, for $j\leqslant d-1$,
$$
(\Cal W-\{\bold p\})\cap\supp\big(\psi_{z_{j}-p_{j}}[-1]\dots\psi_{z_0-p_0}[-1]\Adot\big)=(\Cal W-\{\bold p\})\cap V(z_0-p_0, \dots, z_j-p_j)\cap\supp\Adot,
$$
and the dimension of these sets is at most $d-j-1$, where, as usual, the dimension of the empty set is taken to be $-\infty$.

 Therefore,
$$
\dim_\bold p\big(V(z_0-p_0, \dots, z_{d-1}-p_{d-1})\cap\supp\Adot\big)\leqslant 0,
$$
and, for all $j$ such
that $0\leqslant j\leqslant n$,
$$
\dim_\bold p(\supp\psi_{z_{j}-p_{j}}[-1]\dots\psi_{z_0-p_0}[-1]\Adot)\leqslant d-j-1;
$$ in particular, if $j\geqslant d$, then $H^\bullet\big(\psi_{z_{j}-p_{j}}[-1]\dots\psi_{z_0-p_0}[-1]\Adot\big)_\bold p=0$. }

\vskip .2in

\noindent{\it Proof}. Using the characterization of $\Adot$-isolating given in 5.10.a, we conclude that there exists an open neighborhood $\Cal W$ of $\bold p$ such that, for all $j\leqslant d-1$, $(\Cal W-\{\bold p\})\cap\supp\big(\Phi^j_{{}_{\Adot, \bold z}}(\bold p)\big)=\emptyset$.  Applying Corollary 3.6 inductively, we arrive at the support equality of the proposition, together with the dimension statement. 

If we let $j=d-1$ in the support equality, we conclude that $$
\dim_\bold p\big(V(z_0-p_0, \dots, z_{d-1}-p_{d-1})\cap\supp\Adot\big)\leqslant 0.
$$

Moreover, the support equality and its accompanying dimension statement immediately imply that,
for all $j$ such
that $0\leqslant j\leqslant d-1$,
$$
\dim_\bold p(\supp\psi_{z_{j}-p_{j}}[-1]\dots\psi_{z_0-p_0}[-1]\Adot)\leqslant d-j-1;
$$
in particular, 
$$
\dim_\bold p(\supp\psi_{z_{d-1}-p_{d-1}}[-1]\dots\psi_{z_0-p_0}[-1]\Adot)\leqslant 0.$$
However, this certainly implies that $\bold p\not\in \supp\psi_{z_{d}-p_{d}}[-1]\dots\psi_{z_0-p_0}[-1]\Adot$, and the remainder of the proposition follows.
\qed

\vskip .3in

The definition of $\Adot$-isolating coordinates uses only the coordinate functions $z_0$ through $z_{d-1}$. On the other hand, Definition 5.1 uses the coordinate $z_d$ in the definition of ${}^k\gamma^d_{{}_{\Adot, \bold z}}(\bold p)$ and also gives  the definition of ${}^k\gamma^j_{{}_{\Adot, \bold z}}(\bold p)$ for $j>d$. However, the following corollary tells us that, if our coordinates are $\Adot$-isolating, then ${}^k\gamma^d_{{}_{\Adot, \bold z}}(\bold p)$ only depends on the coordinates $(z_0, \dots, z_{d-1})$ and that, for $j>d$, ${}^k\gamma^j_{{}_{\Adot, \bold z}}(\bold p)$ is trivial.

\vskip .3in

\noindent{\bf Corollary 5.15}.  {\it If the coordinates $\bold z$ are  $\Adot$-isolating at $\bold p$, then $$ {}^k\gamma^d_{{}_{\Adot, \bold z}}(\bold p) \ \cong \
H^k\big(\psi_{z_{d-1}-p_{d-1}}[-1]\psi_{z_{d-2}-p_{d-2}}[-1]\dots\psi_{z_0-p_0}[-1]\Adot\big)_\bold p;
$$

and, for $j>d$, $ {}^k\gamma^j_{{}_{\Adot, \bold z}}(\bold p) \ := \ 0$.
}

\vskip .2in

\noindent{\it Proof}. 
Let $\Adot := \psi_{z_{d-1}-p_{d-1}}[-1]\psi_{z_{d-2}-p_{d-2}}[-1]\dots\psi_{z_0-p_0}[-1]\Adot$. By Proposition 5.14, $\bold p$
is either an isolated point of the support of $\Adot$ or is not in the support at all. In either case, the nearby cycles
$\psi_{z_d-p_d}[-1]\Adot$ are zero in a neighborhood of $\bold p$; therefore, ${}^k\gamma^j_{{}_{\Adot, \bold
z}}(\bold p)=0$ for $j> d$. 

As for
${}^k\gamma^d_{{}_{\Adot, \bold z}}(\bold p)$, consider the fundamental distinguished triangle
$$
\Adot_{|_{V(z_d-p_d)}}[-1] \ \rightarrow \ \psi_{z_d-p_d}[-1]\Adot \ \rightarrow \ \phi_{z_d-p_d}[-1]\Adot\
@>[1]>>\Adot_{|_{V(z_d-p_d)}}[-1].
$$ As $\psi_{z_d-p_d}[-1]\Adot = 0$ near $\bold p$, it follows immediately that $H^k(\Adot)_\bold p\cong
H^k(\phi_{z_d-p_d}[-1]\Adot)_\bold p$, i.e., that
$$H^k\big(\psi_{z_{d-1}-p_{d-1}}[-1]\psi_{z_{d-2}-p_{d-2}}[-1]\dots\psi_{z_0-p_0}[-1]\Adot\big)_\bold p= {}^k\gamma^d_{{}_{\Adot, \bold z}}(\bold p) .\qed
$$ 

\vskip .3in

\noindent{\it Remark 5.16}. As we mentioned before, being $\Adot$-isolating is a condition on only the first $d$ coordinates. Corollary 5.15 tells us that when the coordinates are $\Adot$-isolating, then the characteristic polar modules also depend on only the first $d$ coordinates.

On the other hand, it is useful to have a characterization of $\Adot$-isolating that does not explicitly use the local dimension $d$ of the support of $\Adot$. Recall from Remark 5.6 that, as germs at $\bold p$, $\Theta^m_{{}_{\Adot, \bold z}} = \supp\Adot$ for all $m\geqslant d$. At the same time, Proposition 5.14 tells us that, if $\bold z$ is $\Adot$-isolating at $\bold p$, then $$
\dim_\bold p\big(V(z_0-p_0, \dots, z_{d-1}-p_{d-1})\cap\supp\Adot\big)\leqslant 0.
$$
From these facts, it follows immediately that, if $\bold z$ is $\Adot$-isolating at $\bold p$, then  $$\dim_\bold p \Big(V(z_0-p_0, \dots, z_{j-1}-p_{j-1})\cap \Theta^j_{{}_{\Adot, \bold z}}\Big)\leqslant 0$$ for all $j$ such that $d\leqslant j\leqslant n$.

What this means is that, while the characterization of $\Adot$-isolating given by 5.10.c is a condition for all $j$ such that $0\leqslant j\leqslant d-1$, it is equivalent to define $\Adot$-isolating by the seemingly stronger condition that, for all $j$ such that $0\leqslant j\leqslant n$, $$\dim_\bold p \Big(V(z_0-p_0, \dots, z_{j-1}-p_{j-1})\cap \Theta^j_{{}_{\Adot, \bold z}}\Big)\leqslant 0.$$

\vskip .3in

\noindent{\bf Corollary 5.17}. {\it 
If $R$ is a domain and $\bold z$ is $\Adot$-isolating, then the Euler characteristic of the stalk of $\Adot$ at $\bold p$ is given by
$$\chi(\Adot)_\bold p = \sum_{0\leqslant j\leqslant d}(-1)^j\big({}^{\operatorname{ord}}\gamma^j_{{}_{\Adot, \bold
z}}(\bold x)\big) = \sum_{j, k}(-1)^{j+k}\operatorname{rk}\big({}^k\gamma^j_{{}_{\Adot, \bold
z}}(\bold x)\big).$$
}

\vskip .2in

\noindent{\it Proof}. This follows immediately from Proposition 5.3, using $m=d$, together with Proposition 5.14.\qed

\vskip .3in

\noindent{\bf Theorem 5.18}. {\it 
If $\bold z$ is  $\Adot$-isolating at $\bold p$ and $\Adot$ is a perverse sheaf on
$Y$, then, for all $j$, ${}^\bullet\gamma^j_{{}_{\Adot, \bold z}}(\bold p)$, is supported only in degree $0$, i.e.,
${}^k\gamma^j_{{}_{\Adot, \bold z}}(\bold p)=0$ if $k\neq 0$.

\vskip .2in

 Moreover, whenever $\bold z$ is  $\Adot$-isolating at $\bold p$ and all of the ${}^\bullet\gamma^j_{{}_{\Adot, \bold z}}(\bold p)$ are supported only in degree $0$, there is a complex of
$R$-modules in which ${}^0\gamma^{-m}_{{}_{\Adot, \bold z}}(\bold p)$ is placed in degree $m$,
$$ 0\rightarrow {}^0\gamma^d_{{}_{\Adot, \bold z}}(\bold p)\rightarrow {}^0\gamma^{d-1}_{{}_{\Adot, \bold z}}(\bold
p)\rightarrow\dots\rightarrow {}^0\gamma^1_{{}_{\Adot, \bold z}}(\bold p)\rightarrow {}^0\gamma^0_{{}_{\Adot, \bold z}}(\bold
p)\rightarrow 0,
$$ whose cohomology is isomorphic to $H^m(\Adot)_\bold p$ in degrees $-d\leqslant m\leqslant 0$; in addition, if 
$m\geqslant 1$ or 
$m\leqslant -d-1$,
then $H^m(\Adot)_\bold p = 0$.

\vskip .2in

More generally, whenever $\bold z$ is  $\Adot$-isolating at $\bold p$, for all $k$, there is a complex of
$R$-modules in which ${}^k\gamma^{-m}_{{}_{\Adot, \bold z}}(\bold p)$ is placed in degree $m$,
$$ 0\rightarrow {}^k\gamma^d_{{}_{\Adot, \bold z}}(\bold p)\rightarrow {}^k\gamma^{d-1}_{{}_{\Adot, \bold z}}(\bold
p)\rightarrow\dots\rightarrow {}^k\gamma^1_{{}_{\Adot, \bold z}}(\bold p)\rightarrow {}^k\gamma^0_{{}_{\Adot, \bold z}}(\bold
p)\rightarrow 0,
$$ whose cohomology is isomorphic to $H^m\big({}^\mu\hskip -.02in H^k(\Adot)\big)_\bold p$ in degrees $-d\leqslant m\leqslant 0$; in addition, if 
$m\geqslant 1$ or 
$m\leqslant -d-1$,
then $H^m\big({}^\mu\hskip -.02in H^k(\Adot)\big)_\bold p = 0$.
}

\vskip .2in

\noindent{\it Proof}.  The first statement follows at once from the facts that the functors $\psi_f[-1]$ and $\phi_f[-1]$ take perverse sheaves to perverse sheaves, and perverse sheaves that are supported at isolated points can have non-zero cohomology only in degree zero.

The last statement follows at once from the first and second, together with Proposition 5.13.

\vskip .2in

We must now prove the second statement. It will be convenient to adopt the convention that, for $j=0$,
$$\psi_{z_{j-1}-p_{j-1}}[-1]\dots\psi_{z_0-p_0}[-1]\Adot=\Adot.$$

\vskip .1in

There are fundamental distinguished triangles given by
$$
\Adot_{|_{V(z_j-p_j)}}[-1] \ \rightarrow \ \psi_{z_j-p_j}[-1]\Adot \ \rightarrow \ \phi_{z_j-p_j}[-1]\Adot\
@>[1]>>\Adot_{|_{V(z_j-p_j)}}[-1].\tag{$\dagger$}
$$

For any $j$ such that
$0\leqslant j\leqslant d$, let
$\Adot_{j-1} :=
\psi_{z_{j-1}-p_{j-1}}[-1]\dots\psi_{z_0-p_0}[-1]\Adot$. 
 The fact that all of the ${}^\bullet\gamma^j_{{}_{\Adot, \bold z}}(\bold p)$ are supported only in degree $0$, combined with
$(\dagger)$ yields that, for
$m\neq 0,1$, 
$$ H^{m-1}(\Adot_{j-1})_\bold p =H^m(\Adot_{j-1}[-1])_\bold p \ \cong H^m(\psi_{z_j-p_j}[-1]\Adot_{j-1})_\bold p \ = \ H^{m}(\Adot_j)_\bold p,\tag{$*$}
$$ and an exact sequence
$$ 0\rightarrow H^{-1}\big(\Adot_{j-1}\big)_\bold p\rightarrow
H^{0}(\Adot_j)_\bold p\rightarrow {}^0\gamma^j_{{}_{\Adot, \bold z}}(\bold p)\rightarrow
H^{0}\big(\Adot_{j-1}\big)_\bold p\rightarrow
H^{1}(\Adot_j)_\bold p\rightarrow 0.\tag{$\ddag$}
$$

An inductive application of $(*)$ implies that if $m\geqslant 1$ or if $m\leqslant -d-1$, then
$$ H^m(\Adot)_\bold p \ \cong \
H^{m+d}\big(\Adot_{d-1}\big)_\bold p,
$$
and Corollary 5.15 tells us that this is isomorphic to
${}^{m+d}\gamma^d_{{}_{\Adot, \bold z}}(\bold p)$.  Now, our hypothesis that ${}^\bullet\gamma^d_{{}_{\Adot, \bold z}}(\bold p)$ is concentrated in
degree $0$ implies that ${}^{m+d}\gamma^d_{{}_{\Adot, \bold z}}(\bold p)=0$; this proves the last part of the proposition.
\vskip .1in

Another inductive application of $(*)$ implies that 
$$H^{1}\big(\Adot_j\big)_\bold p \ \cong \ {}^{d-j}\gamma^d_{{}_{\Adot, \bold
z}}(\bold p) = 0,
$$ where the last equality follows from our hypotheses if $d\neq j$ , and if $d=j$, $H^{1}\big(\Adot_j\big)_\bold p = 0$ by the last part of Proposition 5.14 . Thus, the exact
sequences given by
$(\ddag)$ reduce to
$$ 0\rightarrow H^{-1}\big(\Adot_{j-1}\big)_\bold p@>a_j>>
H^{0}(\Adot_j)_\bold p@>b_j>>
 {}^0\gamma^j_{{}_{\Adot, \bold z}}(\bold p)@>c_j>>
H^{0}\big(\Adot_{j-1}\big)_\bold p\rightarrow 0.
$$

Therefore, it is immediate that there is a complex given by
$$ 0\rightarrow {}^0\gamma^d_{{}_{\Adot, \bold z}}(\bold p)@>b_{d-1}\circ c_d>> {}^0\gamma^{d-1}_{{}_{\Adot, \bold
z}}(\bold p)@>b_{d-2}\circ c_{d-1}>>\dots@>b_1\circ c_2>> {}^0\gamma^1_{{}_{\Adot, \bold z}}(\bold
p)@>b_0\circ c_1>> {}^0\gamma^0_{{}_{\Adot, \bold z}}(\bold p)\rightarrow 0,
$$ where ${}^0\gamma^{-m}_{{}_{\Adot, \bold z}}(\bold p)$ stands in degree $m$, and the cohomology of this complex in degree
$m\leqslant -1$ is given by
$$
\frac{\operatorname{ker}(b_{m-1}\circ c_m)}{\operatorname{im}(b_{m}\circ c_{m+1})}=
\frac{\operatorname{ker}(b_{m-1}\circ c_m)}{\operatorname{im}(b_{m})} =
\frac{\operatorname{ker}(b_{m-1}\circ c_m)}{\operatorname{ker}(c_{m})}\cong \operatorname{ker}(b_{m-1})=
\operatorname{im}(a_{m-1})\cong  H^{-1}\big(\Adot_{m-2}\big)_\bold p.$$ 
By another inductive application of $(*)$, this last module is isomorphic to $H^{-m}(\Adot)_\bold p$. The remaining case, in
which $m=0$, is trivial.\qed

\vskip .3in

\noindent{\bf Definition 5.19}.  If $\bold z$ is $\Adot$-isolating at $\bold p$, then we refer to the second chain complex of Theorem 5.18 as the {\it degree $k$ Zawatsky complex (or the degree $k$ Z-complex) of $\Adot$ at $\bold p$
with respect to
$\bold z$}.

\vskip .3in

\noindent{\bf Definition 5.20}. If $\bold z$ is $\Adot$-isolating on the open set $\Cal W$, then, in light of 5.10.d, we may define the {\it graded, enriched, $j$-dimensional characteristic polar cycle of $\Adot$ with respect to $\bold z$, ${}^\bullet\Gamma^j_{{}_{\Adot, \bold z}}$, inside $\Cal W$}  by
$$
{}^k\Gamma^j_{{}_{\Adot, \bold z}}:= \nu_*\Big(\Bbb P(\gecc^k(\Adot)) \ \odot\ \Cal W\times\Bbb
P^j\times\{\bold 0\}\Big).
$$

 If $\bold z$ is $\Adot$-isolating at $\bold p$, then there exists an open neighborhood $\Cal W$ of $\bold p$ on which $\bold z$ is $\Adot$-isolating; thus, the {\it germ of ${}^k\Gamma^j_{{}_{\Adot, \bold z}}$ at $\bold p$} is well-defined.

\vskip .3in

\noindent{\it Remark 5.21}. Note that, if $\bold z$ is $\Adot$-isolating at $\bold p$, then, for all $j$, there is an equality of germs of sets: $\Gamma^j_{{}_{\Adot, \bold z}} =|{}^\bullet\Gamma^j_{{}_{\Adot, \bold z}}|$.

\vskip .3in

\noindent{\bf Lemma 5.22}. {\it Suppose that $\bold z$ is $\Adot$-isolating at $\bold p$. Let $\tilde{\bold z}$ denote the rotated coordinate system $(z_1, z_2, \dots, z_n, z_0)$.  
Then, $\tilde{\bold z}$ is $\big(\psi_{z_0-p_0}[-1]\Adot\big)$-isolating at $\bold p$, and
for all $j$ such that $1\leqslant j\leqslant n$, for all $k$, there is an equality of germs of graded enriched cycles given by
$$
{}^k\Gamma^{j-1}_{{}_{\psi_{z_0-p_0}[-1]\Adot, \tilde{\bold z}}} \ = \ {}^k\Gamma^j_{{}_{\Adot, \bold z}}\odot V(z_0-p_0).
$$
}

\vskip .1in

\noindent{\it Proof}. That $\tilde{\bold z}$ is $\big(\psi_{z_0-p_0}[-1]\Adot\big)$-isolating at $\bold p$ is immediate from the definition, once we show that $\dim_\bold p\supp(\psi_{z_0-p_0}[-1]\Adot)\leqslant -1+\dim_\bold p\supp\Adot$. 

If $\bold p\not\in\supp\Adot$, then $\dim_\bold p\supp(\psi_{z_0-p_0}[-1]\Adot)=\dim_\bold p\supp\Adot=-\infty$. If $\dim_\bold p\supp\Adot = 0$, then $\bold p\not\in\supp(\psi_{z_0-p_0}[-1]\Adot)$, i.e., $\dim_\bold p\supp(\psi_{z_0-p_0}[-1]\Adot)=-\infty$. Therefore, throughout the remainder of the proof, we suppose that $\dim_\bold p\supp\Adot \geqslant 1$.

Since $\bold z$ is $\Adot$-isolating at $\bold p$, $\dim_\bold p\supp(\phi_{z_0-p_0}[-1]\Adot)\leqslant 0$, and so there exists an open neighborhood $\Cal W$ of $\bold p$ such that $\phi_{z_0-p_0}[-1]\Adot$ is zero when restricted to $\Cal W-\{\bold p\}$. That $$\dim_\bold p\supp(\psi_{z_0-p_0}[-1]\Adot)\leqslant -1+\dim_\bold p\supp\Adot$$ now follows at once from Corollary 3.6.

\vskip .2in

Throughout the remainder of the proof, we shall work in the neighborhood $\Cal W$ from above.

\vskip .1in

We have equalities of enriched cycles:
$$
{}^k\Gamma^{j-1}_{{}_{\psi_{z_0-p_0}[-1]\Adot, \tilde{\bold z}}} \ = \ \nu_*\big(\Bbb P(\gecc^k(\psi_{z_0-p_0}[-1]\Adot))\ \odot\ (\Cal W\times\{0\}\times\Bbb P^{j-1}\times\{\bold 0\})\big),
$$
which, by Theorem 3.2, is equal to
$$
\nu_*\Big(\Bbb P\big((T^*_{{}_{z_0, \Adot}}\Cal U)^k)\big)\ \odot\  V(z_0-p_0)\ \odot\ (\Cal W\times\{0\}\times\Bbb P^{j-1}\times\{\bold 0\})\Big)=
$$
$$
V(z_0-p_0)\ \odot\ \nu_*\Big(\Bbb P\big((T^*_{{}_{z_0, \Adot}}\Cal U)^k)\big)\ \odot\ (\Cal W\times\{0\}\times\Bbb P^{j-1}\times\{\bold 0\})\Big),
$$
where the last equality follows from the projection formula. Thus, we would be finished if we could  show that
$$
\nu_*\Big(\Bbb P\big((T^*_{{}_{z_0, \Adot}}\Cal U)^k)\big)\ \odot\ V(w_0)\odot (\Cal W\times\Bbb P^{j}\times\{\bold 0\})\Big) \ = \ 
\nu_*\Big(\Bbb P(\gecc^k(\Adot)) \ \odot\ (\Cal W\times\Bbb
P^j\times\{\bold 0\})\Big).\tag{$\dagger$}
$$
As $j\geqslant 1$, both sides of $(\dagger)$ yield sets whose dimension at $\bold p$ is at least $1$. Therefore, it suffices to prove $(\dagger)$ holds on $\Cal W^\prime:=\Cal W-\{\bold p\}$.

Suppose that $\gecc^k(\Adot)= \sum_\beta E_\beta\big[\overline{T^*_{R_\beta}\Cal U}\big]$. As $\dim_\bold p\phi_{z_0-p_0}[-1]\Adot\leqslant 0$, we may apply Lemma 5.4 to conclude that, above $\Cal W^\prime$, 
$$\Bbb P\big((T^*_{{}_{z_0, \Adot}}\Cal U)^k)\big) = \sum_\beta E_\beta\,\big[\Bbb P\big(\overline{T^*_{R_\beta}\Cal U}+<dz_0>\big)\big].$$
Thus, to prove $(\dagger)$, it suffices to show that, for $\Adot$-visible strata $R_\beta$,
$$
\nu_*\Big(\Bbb P\big(\overline{T^*_{R_\beta}\Cal U}+<dz_0>\big)\big)\ \cdot\ V(w_0)\cdot (\Cal W^\prime\times\Bbb P^{j}\times\{\bold 0\})\Big) \ = \ 
\nu_*\Big(\Bbb P(\overline{T^*_{R_\beta}\Cal U}) \ \odot\ (\Cal W^\prime\times\Bbb
P^j\times\{\bold 0\})\Big).\tag{$\ddagger$}
$$
However, this is easy.

Let $\tau:\Cal W^\prime\times(\Bbb P^n-\{[1:\bold 0]\})\rightarrow \Cal W^\prime\times\{0\}\times\Bbb P^{n-1}$ denote the projection, and note that, over $\Cal W^\prime$, Theorem 3.4 implies that $\Bbb P(\overline{T^*_{R_\beta}\Cal U})\subseteq\Cal W^\prime\times(\Bbb P^n-\{[1:\bold 0]\})$. Now, one notes that $\tau_*\big(\Bbb P(\overline{T^*_{R_\beta}\Cal U})\big)$ is precisely equal to the transverse intersection $\Bbb P\big(\overline{T^*_{R_\beta}\Cal U}+<dz_0>\big)\big)\ \cdot\ V(w_0)$, and then $(\ddagger)$ follows from the projection formula.\qed

\vskip .4in

Finally, we can prove the {\bf Fundamental Theorem of Characteristic Polar Modules}:

\vskip .3in

\noindent{\bf Theorem 5.23}. {\it 
If  $\bold z$ is $\Adot$-isolating at $\bold p$ , then for all $j$ and $k$, where $0\leqslant j\leqslant d$, ${}^k\Gamma^j_{{}_{\Adot, \bold z}}$ properly intersects $V(z_0-p_0, \dots, z_{j-1}-p_{j-1})$ at $\bold p$ and
$$ {}^k\gamma^j_{{}_{\Adot, \bold z}}(\bold p) \ \cong \ \big({}^k\Gamma^j_{{}_{\Adot, \bold z}} \ \odot \ V(z_0-p_0, \dots, z_{j-1}-p_{j-1})\big)_\bold p,$$ 
where, when $j=0$, we mean that ${}^k\gamma^0_{{}_{\Adot, \bold z}}(\bold p) \ = \ \big({}^k\Gamma^0_{{}_{\Adot, \bold z}} \big)_\bold p$.}

\vskip .2in

\noindent{\it Proof}. The proper intersection statement follows immediately from Theorem 5.10.d.  The remainder of the proof is by induction on $j$.

\vskip .1in

When $j=0$, the claim reduces to $$H^k(\phi_{z_0-p_0}[-1]\Adot)_\bold p\cong \Big(\nu_*\big(\Bbb P(\gecc^k(\Adot)) \ \odot\ \Cal W\times\Bbb
P^0\times\{\bold 0\}\big)\Big)_\bold p.$$
This follows immediately from Theorem 3.3.

\vskip .1in

Now, suppose that $j\geqslant 1$ and that the theorem holds for $j-1$ (for arbitrary $\Adot$ and $\bold z$). Then, by Proposition 5.2, our inductive hypothesis, and Lemma 5.22 (in that order),
$$
 {}^k\gamma^j_{{}_{\Adot, \bold z}}(\bold p) =  {}^k\gamma^{j-1}_{{}_{\psi_{z_0-x_0}[-1]\Adot, \tilde{\bold z}}}(\bold p) \ \cong \ \big({}^k\Gamma^{j-1}_{{}_{\Adot, \tilde{\bold z}}} \ \odot \ V(z_1-p_1, \dots, z_{j-1}-p_{j-1})\big)_\bold p \ \cong$$ 
$$
\big({}^k\Gamma^j_{{}_{\Adot, \bold z}} \ \odot \ V(z_0-p_0) \ \odot \ V(z_1-p_1, \dots, z_{j-1}-p_{j-1})\big)_\bold p \ \cong
$$ 
$$
\big({}^k\Gamma^j_{{}_{\Adot, \bold z}} \ \odot \ V(z_0-p_0, \dots, z_{j-1}-p_{j-1})\big)_\bold p.\qed
$$ 

\vskip .3in

\noindent{\it Remark 5.24}. It is useful at this point to step back from all of our definitions and notation, and to describe the result of Theorem 5.23 in terms of ordinary intersection theory. What we have proved throughout this section is the following.

\vskip .1in

The condition that, for all $j<\dim_\bold p\big(\supp\Adot\big)$, $$\dim_\bold p\supp\Big(\phi_{z_j-p_j}[-1]\psi_{z_{j-1}-p_{j-1}}[-1]\dots\psi_{z_0-p_0}[-1]\Adot\Big)\leqslant 0$$  is equivalent to: there exists an open neighborhood $\Cal W$ of $\bold p$ in $\Cal U$ such that, for all strata $R_\beta$ which have non-trivial normal modules (i.e., $\Bbb H^\bullet(\Bbb N_\beta, \Bbb L_\beta; \ \Adot)\neq0$), for all $j\leqslant n$, $\Bbb P(\overline{T^*_{R_\beta}\Cal U})$ properly intersects $\Cal W\times\Bbb P^j\times\{\bold 0\}$ in $\Cal U\times\Bbb P^n$, and 
$$
\dim_\bold p \left\{V(z_0-p_0, \dots, z_{j-1}-p_{j-1})\ \cap\ \nu_*\Big(\Bbb P(\overline{T^*_{R_\beta}\Cal U})\ \cdot \ \big(\Cal W\times\Bbb P^j\times\{\bold 0\}\big)\Big)\right\}\leqslant 0.
$$
Moreover, whenever these equivalent hold, for all $j\leqslant n$, for all $k$,
$$H^k\big(\phi_{z_j-p_j}[-1]\psi_{z_{j-1}-p_{j-1}}[-1]\dots\psi_{z_0-p_0}[-1]\Adot\big)_\bold p \ \cong \ \bigoplus \left(H^{k-d_\beta}(\Bbb N_\beta, \Bbb L_\beta; \ \Adot)\right)^{m_\beta},\tag{$\dagger$}$$
where $d_\beta:=\dim R_\beta$, the sum is over those $\beta$ such that $d_\beta\geqslant j$ and $R_\beta$ has non-trivial normal modules, and 
$$m_\beta:=\left(V(z_0-p_0, \dots, z_{j-1}-p_{j-1})\ \cdot\ \nu_*\Big(\Bbb P(\overline{T^*_{R_\beta}\Cal U})\ \cdot \ \big(\Cal W\times\Bbb P^j\times\{\bold 0\}\big)\Big)\right)_\bold p.$$

\vskip .1in

If the coordinates $\bold z$ are generic enough at $\bold p$, and still assuming that $d_\beta\geqslant j$, then $m_\beta$ is equal to the multiplicity of the $j$-dimensional absolute polar variety of $\overline{R_\beta}$ at $\bold p$. Note, however, while $\Adot$-isolating coordinates at $\bold p$ are $\Adot$-isolating throughout a neighborhood of $\bold p$, that -- unless $\bold p$ is a non-singular point of $\overline{R_\beta}$ -- one may {\bf not} choose one set of coordinates which are generic enough to yield the polar multiplicities of $\overline{R_\beta}$ throughout a neighborhood of $\bold p$.

One interesting observation that one can make from looking at $(\dagger)$ is that, if all of the normal modules of strata are free (or merely torsion-free), then so are all of the characteristic polar modules at each point.

\vskip .3in

\noindent\S6. {\bf L\^e-Vogel Cycles and Numbers}

\vskip .2in

We wish to apply the results of the previous section to the case where $Y:=V(f)$ and $\Adot:=\phi_f[-1]\Fdot$. Corollary 5.17 and Theorem 5.18 tell us that the characteristic polar modules of $\phi_f[-1]\Fdot$ provide a great deal of information about the stalk cohomology of $\phi_f[-1]\Fdot$. Thus, what we need to do is to find an algebraic method for calculating $\nu_*\Big(\Bbb P(\gecc^k(\phi_f[-1]\Fdot)) \ \odot\ \Cal U\times\Bbb
P^j\times\{\bold 0\}\Big)$. We also need to investigate when coordinates are $\phi_f[-1]\Fdot$-isolating.

Note that, if $f\equiv 0$, then $\phi_f[-1]\Fdot\cong \Fdot$, and therefore the relative results of this section also apply to the absolute situation of the previous section.

\vskip .2in

We recall the following definition from [{\bf M1}].

\vskip .2in

\noindent{\bf Definition 6.1}. The {\it $\Fdot$-critical locus of $f$}, $\Sigma_{{}_{\Fdot}}f$, is equal to $\{\bold x\in X\ |\ H^\bullet(\phi_{f-f(\bold x)}[-1]\Fdot)_\bold x\neq 0\}$. Hence, the closure, $\dsize\overline{\Sigma_{{}_{\Fdot}}f}=\bigcup_v\supp(\phi_{f-v}[-1]\Fdot)$.

\vskip .3in

\noindent{\bf Proposition 6.2}. {\it For all $v\in\Bbb C$, 
$$
V(f-v)\cap \Sigma_{{}_{\Fdot}}f \ = \ \eta(|\gecc^\bullet(\phi_{f-v}[-1]\Fdot)|) \ = \ \nu(|\Bbb P\big(\gecc^\bullet(\phi_{f-v}[-1]\Fdot)\big)|).
$$
}

\vskip .1in

\noindent{\it Proof}. This follows immediately from the third equality of Proposition 2.5, applied to the complex $\phi_{f-v}[-1]\Fdot$.
\qed

\vskip .3in

Note that, with our new notation, $d=\dim_\bold p\overline{\Sigma_{{}_{\Fdot}}f}$.

\vskip .3in

At long last, we can define our generalizations of the L\^e cycles and L\^e numbers to analytic functions with arbitrarily singular domains.

\vskip .3in

\noindent{\bf Definition 6.3}. For all $j$ such that $0\leqslant j\leqslant n$, for all $k$, we let
$$
{}^k\lambda^j_{{}_{f, \bold z}}(\bold p;\ \Fdot) \ := \ {}^k\gamma^j_{{}_{\phi_f[-1]\Fdot, \bold z}}(\bold p);
$$
that is
$$
{}^k\lambda^j_{{}_{f, \bold z}}(\bold p;\ \Fdot) \ := \ H^k\big(\phi_{z_j-p_j}[-1]\psi_{z_{j-1}-p_{j-1}}[-1]\dots\psi_{z_0-p_0}[-1]\phi_f[-1]\Fdot\big)_\bold p,
$$
and we refer to ${}^k\lambda^j_{{}_{f, \bold z}}(\bold p;\ \Fdot)$ as the {\it degree $k$, $j$-dimensional L\^e-Vogel (L\^eVo) module of $f$ at $\bold p$ with respect to $\bold z$ with coefficients in $\Fdot$}.

\vskip .1in

 For all $j$ such that  $0\leqslant j\leqslant n$, we
define {\it $\Omega^j_{{}_{f, \bold z}}(\Fdot)$ } by 
$$
 \Omega^j_{{}_{f, \bold z}}(\Fdot) \ := \  \nu\Big(|\Bbb P(\gecc^\bullet(\phi_f[-1]\Fdot))| \ \cap \ (\Cal U\times\Bbb
P^j\times\{\bold 0\})\Big),
$$
and we define the {\it $j$-dimensional L\^e-Vogel (L\^e-Vogel) set of $f$ with respect to $\bold z$ with coefficients in $\Fdot$}, $\Lambda^j_{{}_{f, \bold z}}(\Fdot)$, to be the union of the $j$-dimensional components of $\Omega^j_{{}_{f, \bold z}}(\Fdot)$.

\vskip .1in

If $\bold z$  are $\phi_f[-1]\Fdot$-isolating at $\bold p$, then we define the germ at $\bold p$ of the graded, enriched cycle ${}^\bullet\Lambda^j_{{}_{f, \bold z}}(\Fdot)$ by
$$
{}^k\Lambda^j_{{}_{f, \bold z}}(\Fdot) \ := \ \ {}^k\Gamma^j_{{}_{\phi_f[-1]\Fdot, \bold z}},
$$
and refer to this as the germ at $\bold p$ of the {\it degree $k$, $j$-dimensional L\^e-Vogel (L\^eVo) cycle of $f$ with respect to $\bold z$ with coefficients in $\Fdot$}.

\vskip .3in

Our earlier results allow us to quickly prove

\vskip .3in

\noindent{\bf Theorem 6.4}. {\it For a generic choice of $\bold z$, the coordinates $\bold z$ are $\phi_f[-1]\Fdot$-isolating at $\bold p$.

\vskip .1in

If  the coordinates $\bold z$ are $\phi_f[-1]\Fdot$-isolating at $\bold p$, then there is an open neighborhood, $\Cal W$, of $\bold p$ on which $\bold z$ are $\phi_f[-1]\Fdot$-isolating.

\vskip .1in

If  the coordinates $\bold z$ are $\phi_f[-1]\Fdot$-isolating at $\bold p$, then, for all $k$, for all $j>\dim_\bold p\supp(\phi_f[-1]\Fdot)$, $ {}^k\lambda^j_{{}_{f, \bold z}}(\bold p;\ \Fdot)=0$.

\vskip .1in

The coordinates  $\bold z$ are $\phi_f[-1]\Fdot$-isolating on an open $\Cal W\subseteq \Cal U$   if and only if, for all $m$ 
such that $0\leqslant m\leqslant n$, $|\Bbb P\big(\gecc^\bullet(\phi_f[-1]\Fdot)\big)|$ properly intersects 
$\Cal W\times \Bbb P^m\times\{\bold 0\}$ in $\Cal U\times\Bbb P^n$ and, for all $\bold x\in\Cal W$, $\dim_\bold x \left(V(z_0-x_0, \dots, z_{m-1}-x_{m-1})\cap\Lambda^m_{{}_{f, \bold z}}(\Fdot)\right)\leqslant 0$.

\vskip .1in

If the coordinates  $\bold z$ are $\phi_f[-1]\Fdot$-isolating on an open $\Cal W\subseteq \Cal U$, then, on $\Cal W$, there is an equality of graded, enriched cycles given by
$$
{}^k\Lambda^j_{{}_{f, \bold z}}(\Fdot)= \nu_*\Big(\Bbb P(\gecc^k(\phi_f[-1]\Fdot)) \ \odot\ \Cal W\times\Bbb
P^j\times\{\bold 0\}\Big),
$$
and, for all $\bold x\in V(f)$, for all $k$ and $j$, and isomorphism of modules 
$$ {}^k\lambda^j_{{}_{f, \bold z}}(\bold x;\ \Fdot) \ \cong \ \big({}^k\Lambda^j_{{}_{f, \bold z}}(\Fdot) \ \odot \ V(z_0-x_0, \dots, z_{j-1}-x_{j-1})\big)_\bold x.$$

\vskip .1in

If  the coordinates $\bold z$ are $\phi_f[-1]\Fdot$-isolating at $\bold p$, then, for all $k$, there is a complex of
$R$-modules in which $ {}^k\lambda^{-m}_{{}_{f, \bold z}}(\bold p;\ \Fdot)$ is placed in degree $m$,
$$ 0\rightarrow{}^k\lambda^{d}_{{}_{f, \bold z}}(\bold p;\ \Fdot)\rightarrow {}^k\lambda^{d-1}_{{}_{f, \bold z}}(\bold p;\ \Fdot)\rightarrow\dots\rightarrow {}^k\lambda^{1}_{{}_{f, \bold z}}(\bold p;\ \Fdot)\rightarrow{}^k\lambda^{0}_{{}_{f, \bold z}}(\bold p;\ \Fdot)\rightarrow 0,
$$ whose cohomology is isomorphic to $H^m\big({}^\mu\hskip -.02in H^k(\phi_f[-1]\Fdot)\big)_\bold p\cong H^m\big(\phi_f[-1]\big({}^\mu\hskip -.02in H^k(\Fdot)\big)\big)_\bold p$ in degrees $-d\leqslant m\leqslant 0$; in addition, if 
$m\geqslant 1$ or 
$m\leqslant -d-1$,
then $$H^m\big({}^\mu\hskip -.02in H^k(\phi_f[-1]\Fdot)\big)_\bold p \cong H^m\big(\phi_f[-1]\big({}^\mu\hskip -.02in H^k(\Fdot)\big)\big)_\bold p= 0.$$

\vskip .1in

If the base ring $R$ is a domain, and the coordinates $\bold z$ are $\phi_f[-1]\Fdot$-isolating at $\bold p$, then
$$
\chi(\phi_f[-1]\Fdot)_\bold p = \sum_{j, k} (-1)^{k+j}\operatorname{rk}\big( {}^k\lambda^j_{{}_{f, \bold z}}(\bold p;\ \Fdot)\big).
$$
}

\vskip .1in

\noindent{\it Proof}. All of the results here follow immediately by applying results of the previous section to our current situation.
\qed

\vskip .3in

The first statement in Theorem 6.4 tells us that $\phi_f[-1]\Fdot$-isolating coordinates are generic, but we would like to have a reasonable idea of just how generic the coordinates need to be. Recall the definition of an $a_{{}_{f, \Fdot}}$ partition from Definition 4.7.

\vskip .3in

\noindent{\bf Theorem 6.5}. {\it  Suppose that $\Cal S$ is an $a_{{}_{f, \Fdot}}$ partition of $X$, and let $\Cal S^\prime:=\{S_\alpha\ |\ S_\alpha\subseteq V(f)\}$. Then, $\Cal S^\prime$ is a $\phi_f[-1]\Fdot$-partition of $V(f)$.

Therefore, if $\bold z$ is essentially transverse to all of the $V(f)$ strata of an $a_{{}_{f, \Fdot}}$ partition, then $\bold z$ is $\phi_f[-1]\Fdot$-isolating.}

\vskip .2in

\noindent{\it Proof}. Suppose that $\bold x\in S_\beta\in\Cal S^\prime$, and that $S_\alpha$ is an $\Fdot$-visible stratum. Then, by Theorem 4.8, there is an inclusion of fibers over $\bold x$ given by
$$
\pi\big({\operatorname{Ex}}_{\imdf}(\con)\big)_\bold x \ \subseteq \ \Bbb P\big(T^*_{{}_{S_\beta}}\Cal U\big)_\bold x.
$$
It follows that 
$$\nu^{-1}(V(f)) \ \cap \ \pi\big({\operatorname{Ex}}_{\imdf}(\con)\big) \ \subseteq \ \bigcup_{S_\beta\in\Cal S^\prime}\Bbb P\big(T^*_{{}_{S_\beta}}\Cal U\big).$$
The proposition is now an immediate consequence of Theorem 3.4.
\qed

\vskip .3in

\noindent{\it Remark 6.6}. Theorem 6.5  explains the importance of {\it good stratifications} and {\it prepolar coordinates} in our earlier work (see, for instance, [{\bf M3}]) on L\^e cycles and numbers. With our current notation and terminology, prepolar coordinates essentially transverse to all of the $V(f)$ strata of an $a_{{}_{f, \Bbb C^\bullet_{{}_\Cal U}}}$ partition

\vskip .3in

The following is the {\bf Fundamental Theorem of L\^e-Vogel Cycles}.

\vskip .2in

\noindent{\bf Theorem 6.7}. {\it Then, the following are equivalent:

\vskip .1in

\noindent a)\hskip .2in  the coordinates $\bold z$ are $\phi_f[-1]\Fdot$-isolating at $\bold p$;

\vskip .1in

\noindent b)\hskip .2in there exists a neighborhood $\Cal W$ of $\bold p$ in $\Cal U$ such that, for all $j$ with $0\leqslant j\leqslant n$, the graded, enriched cycle $\pi_*\big({\operatorname{Ex}}_{\operatorname{im}d\tilde f}(\gecc^\bullet(\Fdot))\big)$ properly intersects $\Cal W\times\Bbb P^j\times\{\bold 0\}$ and 
$$\dim_\bold p V(z_0-p_0, \dots, z_{j-1}-p_{j-1})\cap \left|\nu_*\Big(\pi_*\big({\operatorname{Ex}}_{\operatorname{im}d\tilde f}(\gecc^\bullet(\Fdot))\big)\odot\Cal W\times\Bbb P^j\times\{\bold 0\}\Big)\right| \ \leqslant \ 0;
$$
\vskip .1in

\noindent c)\hskip .2in there exists a neighborhood $\Cal W$ of $\bold p$ in $\Cal U$ such that, for all $j$ with $0\leqslant j\leqslant n$, the graded, enriched cycle ${\operatorname{Ex}}_{\operatorname{im}d\tilde f}(\gecc^k(\Fdot))$ properly intersects $\Cal W\times\Bbb C^{n+1}\times\Bbb P^j\times\{\bold 0\}$ and 
$$\dim_\bold p V(z_0-p_0, \dots, z_{j-1}-p_{j-1})\cap \left|\eta_*\Big(\tau_*\Big({\operatorname{Ex}}_{\operatorname{im}d\tilde f}(\gecc^k(\Fdot)) \ \odot \ \Cal U\times\Bbb C^{n+1}\times\Bbb P^j\times\{\bold 0\}\Big)\Big)\right| \ \leqslant \ 0.$$

\vskip .2in

 Whenever the equivalent conditions above hold, there are equalities of germs at $\bold p$ of enriched cycles given by
$$
{}^k\Lambda^j_{{}_{f, \bold z}}(\Fdot) \ = \ \nu_*\Big(\pi_*\big({\operatorname{Ex}}_{\operatorname{im}d\tilde f}(\gecc^k(\Fdot))\big) \ \odot \ \Cal U\times\Bbb P^j\times\{\bold 0\}\Big)=
$$
$$
\eta_*\Big(\tau_*\Big({\operatorname{Ex}}_{\operatorname{im}d\tilde f}(\gecc^k(\Fdot)) \ \odot \ \Cal U\times\Bbb C^{n+1}\times\Bbb P^j\times\{\bold 0\}\Big)\Big),
$$
where, in the last expression, $\eta_*$ induces an isomorphism onto its image, i.e., if we let
$$
\sum_V M_V[V] \ = \ \tau_*\Big({\operatorname{Ex}}_{\operatorname{im}d\tilde f}(\gecc^k(\Fdot)) \ \odot \ \Cal U\times\Bbb C^{n+1}\times\Bbb P^j\times\{\bold 0\}\Big),
$$
then, for all $V$ for which $M_V\neq 0$, $\eta_{|_V}$ is an isomorphism onto its image, and $\eta_*(\sum_V M_V[V])= \sum_V M_V[\eta(V)]$.
}

\vskip .2in

\noindent{\it Proof}. The equivalence of a) and b) follows immediately from 5.10. The equivalence of b) and c) follows trivially from the fact that $\nu\circ\pi=\eta\circ\tau$.

\vskip .2in

As for the formula, the first equality is immediate from Theorems 6.4 and 3.4. 

\vskip .1in

Over $\imdf\times\Bbb P^n$, $\pi$ has an inverse  $\xi:\imdf\times\Bbb P^n\rightarrow\Cal U\times\Bbb C^{n+1}\times\Bbb P^n$ given by $\xi(\bold x, [\bold u])=(\bold x, d_\bold x\tilde f, [\bold u])$. Therefore, elementary intersection theory implies that 
$$
\pi_*\big({\operatorname{Ex}}_{\operatorname{im}d\tilde f}(\gecc^k(\Fdot))\big) \ \odot \ \Cal U\times\Bbb P^j\times\{\bold 0\} \ = \ \pi_*\Big({\operatorname{Ex}}_{\operatorname{im}d\tilde f}(\gecc^k(\Fdot)) \ \odot \ \Cal U\times\Bbb C^{n+1}\times\Bbb P^j\times\{\bold 0\}\Big).
$$ Now, the second equality follows from the fact that $\nu\circ\pi=\eta\circ\tau$.

\vskip .1in

The final isomorphism claim follows from the fact that $$\tau\Big({\operatorname{Ex}}_{\operatorname{im}d\tilde f}(\gecc^k(\Fdot)) \ \odot \ \Cal U\times\Bbb C^{n+1}\times\Bbb P^j\times\{\bold 0\}\Big)$$ lies inside $\imdf$, and -- as with $\pi$ in the above paragraph -- $\eta$ is invertible over $\imdf$.
\qed

\vskip .4in

The importance of Theorem 6.7 is that the general theory of Vogel cycles (see [{\bf G1}], [{\bf G2}], [{\bf V}], and [{\bf M7}]) gives one an algorithm for calculating $\tau_*\Big({\operatorname{Ex}}_{\operatorname{im}d\tilde f}(\gecc^k(\Fdot)) \ \odot \ \Cal U\times\Bbb C^{n+1}\times\Bbb P^j\times\{\bold 0\}\Big)$. This is the algorithm that we stated in a slightly more general form in Section 3 as Theorem 3.8. The corollary below follows immediately from 3.8; hence, we present this corollary without proof.

\vskip .3in

\noindent{\bf Corollary 6.8}. {\it  Suppose that $\bold z$ are $\phi_f[-1]\Fdot$-isolating at  $\bold p$. Then, with respect to $\bold z$, the L\^eVo modules of $f$ at $\bold p$ with coefficients in $\Fdot$, and the germs at $\bold p$ of the corresponding L\^eVo  cycles may be calculated via the following process:

\vskip .1in

We assume that we are working over a sufficiently small neighborhood of $\bold p$ so that the coordinates are $\phi_f[-1]\Fdot$-isolating on the entire neighborhood.

\vskip .1in

Let $\Pi^{n+1}:= \gecc^k(\Fdot)$. Then, $\Pi^{n+1}$ properly intersects $\dsize V\left(w_n-\frac{\partial\tilde f}{\partial z_n}\right)$, and we may consider the enriched cycle defined by the intersection $$\sum_V M_V[V] \ := \ \Pi^{n+1}\odot V\left(w_n-\frac{\partial\tilde f}{\partial z_n}\right);$$ this enriched cycle may have some components contained in $\imdf$ and some components not contained in $\imdf$. Let $\dsize\Pi^n:= \sum_{V\not\subseteq\imdf} M_V[V]$ and let $\dsize\Delta^n:=\sum_{V\subseteq\imdf}M_V[V]$.

\vskip .1in

Now, proceed inductively: if we have $\Pi^{j+1}$, then $\dsize V\left(w_j-\frac{\partial\tilde f}{\partial z_j}\right)$ properly intersects $\Pi^{j+1}$, and we define $\Pi^j$ and $\Delta^j$ by the equality
$$
\Pi^{j+1}\odot V\left(w_j-\frac{\partial\tilde f}{\partial z_j}\right)  \ = \ \Pi^j+\Delta^j,$$
where no component of $\Pi^j$ is contained in $\imdf$, and every component of $\Delta^j$ is contained in $\imdf$.

\vskip .1in

Continue with this process until one obtains $\Pi^0$ and $\Delta^0$.

\vskip .2in

Then, for all $j$, as germs at $\bold p$, ${}^k\Lambda^j_{{}_{f, \bold z}}(\Fdot) = \eta_*(\Delta^j)$ and  
$${}^k\lambda^j_{{}_{f, \bold z}}(\bold p;\ \Fdot) = V(z_0-p_0, \dots, z_{j-1}-p_{j-1})\odot \eta_*(\Delta^j).
$$
}

\vskip .1in

\noindent{\it Remark 6.9}.  It would, of course, be nice if we could use the process of Corollary 6.8 to decide whether or not the coordinates are $\phi_f[-1]\Fdot$-isolating. That is, one might hope that, if all of the intersections in Corollary 6.8 are proper -- including the intersections of $V(z_0-p_0, \dots, z_{j-1}-p_{j-1})$ and  $\eta_*(\Delta^j)$ -- then the coordinates are $\phi_f[-1]\Fdot$-isolating at $\bold p$. While we can not prove this in general, or find a counterexample, we do prove it below in the cases where $d\leqslant 2$. 

In the general case, we must use the more unmanageable condition of requiring the coordinates to be essentially transverse to all of the $V(f)$ strata of an $a_{{}_{f, \Fdot}}$ partition (recall Theorem 6.5).

\vskip .2in

First, we need an easy transversality lemma. For notational ease, we concentrate our attention at the origin.

\vskip .2in

\noindent{\bf Lemma 6.10}. {\it Suppose that $Z$ is an analytic subset of $\Cal U$, and that $\dim_\bold 0 V(z_0, \dots, z_j)\cap Z\leqslant 0$.

Then, for generic $\eta\in \Bbb P^j\times\{\bold 0\}$, there exists a open neighborhood $\Cal W$ of $\bold 0$ in $\Cal U$ such that
$$
\Bbb P\big(\overline{T^*_{{}_{Z_{\operatorname{reg}}}}\Cal U}\big)\cap (\Cal W\times\{\eta\})\subseteq \{\bold 0\}\times\Bbb P^n.
$$
}

\vskip .1in

\noindent{\it Proof}. If $\bold  0\not\in Z$, then the result is trivial. So, assume that  $\dim_\bold 0 V(z_0, \dots, z_j)\cap Z=0$. Let $\xi:\Bbb C^{n+1}\rightarrow \Bbb C^{j+1}$ denote that projection onto the first $j+1$ coordinates.

As $\dim_\bold 0 V(z_0, \dots, z_j)\cap Z=0$, there exists an open neighborhood $\Cal V$ of $\bold 0$ in $\Cal U$, and an open neighborhood $\Cal Q$ of $\bold 0$ in $\Bbb C^{j+1}$ such that the restriction of $\xi$, $\tilde \xi: \Cal V\cap Z\rightarrow\Cal Q$, is a finite map. As $\tilde \xi$ is proper, we may Whitney stratify $\Cal V\cap Z$ and $\operatorname{im}\tilde \xi$ in such a way that $\tilde \xi$ becomes a stratified map (see I.1.7 of [{\bf G-M2}]).

Now, a generic hyperplane through the origin in $\Bbb C^{j+1}$ will transversely intersect all of the strata of $\operatorname{im}\tilde \xi$ near the origin, except possibly at the origin itself. It follows that, for generic $[a_0:\dots:a_j]\in \Bbb P^j$, $V(a_0z_0+\dots+a_jz_j)$ transversely intersects all of the strata of $\Cal V\cap Z$ near the origin, except perhaps at the origin itself. Fix such an $[a_0:\dots:a_j]$.

We proceed by contradiction. Suppose there were a parameterized analytic curve $\bold q(t)\in Z$ such that $\bold q(0)=\bold 0$, $\bold q(t)\neq 0$ for $t\neq 0$, and $$[a_0dz_0+\dots +a_jdz_j]\in\Big(\Bbb P\big(\overline{T^*_{{}_{Z_{\operatorname{reg}}}}\Cal U}\big)\Big)_{\bold q(t)}.$$
For $t$ sufficiently small and unequal to $0$, $\bold q(t)$ must be contained in one of the Whitney strata, $M$, of $Z$. By Whitney's condition a), we have that, for $t$ small and unequal to $0$, 
$$[a_0dz_0+\dots +a_jdz_j]\in\big(\Bbb P\big(T^*_{{}_{M}}\Cal U\big)\big)_{\bold q(t)}.\tag{$\dagger$}$$
As $\bold q^\prime(t)\in T_{\bold q(t)}M$, $(\dagger)$ implies that $a_0q_0^\prime(t) + \dots + a_jq_j^\prime(t)\equiv 0$ for $t\neq 0$. As $\bold q(0)=\bold 0$, it follows that $a_0q_0(t) + \dots + a_jq_j(t)\equiv 0$, i.e., that $\bold q(t)\in V(a_0z_0+\dots+a_jz_j)$. However, this means that $(\dagger)$ contradicts the fact that $V(a_0z_0+\dots+a_jz_j)$ transversely intersects $M$ near the origin.
\qed

\vskip .3in

The proof of the following proposition is not particularly difficult, but it does break up into a number of cases.

\vskip .3in

\noindent{\bf Proposition 6.11}. {\it Assume that we are in the setting of Corollary 6.9, except that we do {\bf not} assume that the  coordinates $\bold z$ are $\phi_f[-1]\Fdot$-isolating at $\bold p$. Assume that $d\leqslant 2$.

If there is a neighborhood of $\bold p$ over which all of the intersections appearing in Corollary 6.9 are proper, including  the intersections $V(z_0-p_0, \dots, z_{j-1}-p_{j-1})\odot \eta_*(\Delta^j)$, then the  coordinates $\bold z$ are, in fact, $\phi_f[-1]\Fdot$-isolating at $\bold p$.
}

\vskip .2in

\noindent{\it Proof}. We use the criterion in 5.10.c for showing that the coordinates $\bold z$ are $\phi_f[-1]\Fdot$-isolating at $\bold p$. For notational convenience, we will assume that $\bold p=\bold 0$. 

We use the notation from 6.8. Note that we have an equality of sets 
$$|\gecc^\bullet(\Fdot)|\cap V\left(w_0-\frac{\partial \tilde f}{\partial z_0}, \dots, w_n-\frac{\partial \tilde f}{\partial z_n}\right) \ = \ \bigcup_j |\Delta^j|;
$$
this is easy to see, or follows rigorously from Proposition I.2.4 and Proposition I.2.15 of [{\bf M7}]. Thus, applying $\eta$ to both sides, and using the first part of Theorem 3.4, we conclude that, near $\bold 0$,
$$
\overline{\Sigma_{{}_{\Fdot}}f} = \bigcup_j \eta(|\Delta^j|).
$$
Therefore, our assumptions on proper intersections imply that 
$$\dim_\bold 0 V(z_0, z_1)\cap\overline{\Sigma_{{}_{\Fdot}}f} \leqslant 0.\tag{$*$}
$$

Let $S_\alpha\in\Cal S$ be an $\Fdot$-visible stratum, and let us simply write $E$ for the exceptional divisor ${\operatorname{Ex}}_{\operatorname{im}d\tilde f}(\overline{T^*_{{}_{S_\alpha}}\Cal U})$. Combining 5.10.c with Theorem 3.4, we see that what we need to show is that: for all $j$,
$$
\dim_\bold 0 V(z_0, \dots, z_{j-1})\cap \left|\nu\Big(\pi(E) \ \cap \ \big(\Cal U\times\Bbb P^j\times\{\bold 0\}\big)\Big)\right|\leqslant 0,\tag{$\dagger$}
$$
where Theorem 3.4 tells us that, as sets, $\pi(E)$ is a union of components of $\bigcup_\beta \Bbb P(\overline{T^*_{{}_{R_\beta}}\Cal U})= \bigcup_\beta \Bbb P(T^*_{{}_{R_\beta}}\Cal U)$
for some Whitney a) stratification, $\{R_\beta\}$, of $V(f)$.
By $(*)$, the condition $(\dagger)$ holds for $j\geqslant 2$. We have only to show that $(\dagger)$ holds for $j=0$ and $j=1$.

\vskip .1in

We proceed by contradiction.

\vskip .1in

\noindent $\bullet$ Case 1: $j=0$.

\vskip .1in

Suppose that there is a parameterized analytic curve $\bold q(t)$ such that $\bold q(0)=\bold 0$, $\bold q(t)\neq \bold 0$ for $t\neq 0$, and, for $t\neq 0$, $$\bold q(t)\in  \left|\nu\Big(\pi(E) \ \cap \ \big(\Cal U\times\Bbb P^0\times\{\bold 0\}\big)\Big)\right|.$$

For $t$ small and unequal to $0$, there exists a single $R_\beta$ such that $\bold q(t)\in R_\beta$. As $\bold q^\prime(t)\in T_{\bold q(t)}R_\beta$ for small $t\neq 0$, it follows that $q_0^\prime(t)\equiv 0$. Since $q_0(0) = 0$, we must have $q_0(t)\equiv 0$, i.e., $\bold q(t)\in V(z_0)$. Therefore, if we can show that $(\dagger)$ holds for $j=1$, the $j=0$ case will follow.

\vskip .1in

\noindent $\bullet$ Case 2: $j=1$.

\vskip .1in

Suppose that $C$ is a component of $\pi(E)\cap (\Cal U\times\Bbb P^1\times\{\bold 0\})$ such that $\dim_\bold 0 V(z_0)\cap \nu(C)\geqslant 1$. Such a component $C$ corresponds to a component $\widetilde C$ of $E\cap (\Cal U\times\Bbb C^{n+1}\times\Bbb P^1\times\{\bold 0\})$, which must be contained in a component $D$ of $\big(\operatorname{Bl}_{\imdf}(\overline{T^*_{{}_{S_\alpha}}\Cal U})\big)\cap (\Cal U\times\Bbb C^{n+1}\times\Bbb P^1\times\{\bold 0\})$. Note that the dimension of $D$ must be at least $2$.

\vskip .2in

If $D\not\subseteq E$, then -- looking at the definition of $\Pi^2$ in 6.8 -- we find that $\tau(D)\subseteq \Pi^2$. Thus, if $D\not\subseteq E$, then 
$$\tau(\widetilde C)\subseteq\tau(D\cap E)\subseteq \Pi^2\cap V\left(\bold w-\frac{\partial \tilde f}{\partial\bold z}\right) = \Delta^1\cup \Delta^0,$$
and so $$V(z_0)\cap\nu(C)= V(z_0)\cap \nu(\pi(\widetilde C))\subseteq V(z_0)\cap\nu(\pi(D\cap E))=$$
$$V(z_0)\cap \eta(\tau(D\cap E)) \subseteq V(z_0)\cap \big(\eta(\Delta^1)\cup \eta(\Delta^0)\big),$$
which, by assumption, has dimension of at most $0$ at the origin. This is a contradiction.

\vskip .2in

It remains only for us to dispose of the case where $D\subseteq E$. Hence, we will be finished if we can show, for every component $C$ of $\pi(E) \ \cap \ \big(\Cal U\times\Bbb P^1\times\{\bold 0\}\big)$ of dimension at least $2$,  that $\dim_\bold 0 V(z_0)\cap\nu(C)\leqslant 0$. 

\vskip .3in

We proceed by contradiction. Suppose that we have an irreducible $C$ such that $\dim C\geqslant 2$ and such that $\dim_\bold 0 V(z_0)\cap\nu(C)\geqslant 1$. 

\vskip .2in

If $\dim_\bold 0 \nu(C) =0$, then there is nothing to show. 

\vskip .2in

Suppose that $\dim_\bold 0 \nu(C) =1$. As $\dim C\geqslant 2$, we must have that $C = \nu(C)\times\Bbb P^1\times\{\bold 0\}$. Let $R_\beta$ be the stratum which contains $\nu(C)-\{\bold 0\}$ in some neighborhood of the origin. Then, over a neighborhood of the origin,
$$\big(\nu(C)-\{\bold 0\}\big)\times\Bbb P^1\times\{\bold 0\} \ \subseteq \  \Bbb P(T^*_{{}_{R_\beta}}\Cal U).$$
However, this means that, along the curve $\nu(C)-\{\bold 0\}$, the tangent space to $R_\beta$ is contained in $\{\bold 0\}\times \Bbb C^{n-1}$. By parameterizing $\nu(C)$, we see that this implies that $\nu(C)\subseteq V(z_0, z_1)$ near the origin. This contradicts $(*)$.

\vskip .2in

Suppose, finally, that $\dim_\bold 0 \nu(C) =2$. Then, there exists a $2$-dimensional stratum $R_\beta$ such that $C\subseteq \Bbb P(\overline{T^*_{{}_{R_\beta}}\Cal U})$, and $\nu(C)=\overline{R_\beta}$.

By $(*)$, $\overline{R_\beta}\not\subseteq V(z_0, z_1)$. Thus, for generic $[a_0: a_1]\in\Bbb P^1$, $V(a_0z_0+a_1z_1)\cap \overline{R_\beta}$ is $1$-dimensional and $\big(V(a_0z_0+a_1z_1)\cap \overline{R_\beta}\big)-\{\bold 0\}\subseteq R_\beta$ . Fix such an $[a_0:a_1]$.  Let $\bold q(t)$ be a parameterization of a component of $V(a_0z_0+a_1z_1)\cap \nu(C)$ such that $\bold q(0)=\bold 0$.  Therefore, for all small $t\neq 0$, $\bold q^\prime(t)\in T_{\bold q(t)}R_\beta$, $a_0q_0^\prime(t)+a_1q_1^\prime(t) = 0$, and $\Bbb P(T^*_{{}_{R_\beta}}\Cal U)_{\bold q(t)}$ contains an element of the form $[b_0(t)dz_0+b_1(t)dz_1]$.

It follows that both $a_0q_0^\prime(t)+a_1q_1^\prime(t) = 0$ and $b_0(t)q_0^\prime(t)+b_1(t)q_1^\prime(t) = 0$ for all small $t\neq 0$; thus, either $q_0^\prime(t)= q_1^\prime(t)=0$ or $[a_0:a_1]=[b_0(t):b_1(t)]$. If $q_0^\prime(t)= q_1^\prime(t)=0$ for an infinite number of $t$ values near $t=0$, then we must have that both $q_0(t)$ and $q_1(t)$ are identically equal to zero; this would contradict $(*)$.

Thus, we must have that, for small $t\neq 0$, $[b_0(t):b_1(t)]$ equals the constant $[a_0:a_1]$, which contradicts the $j=1$ case of Lemma 6.10.\qed

\vskip .3in

\noindent\S7. {\bf An Example}

\vskip .2in

We will now use Proposition 6.11 and Corollary 6.8 to perform calculations in an example in which $\dim\Sigma_{{}_{\Fdot}}f=1$.

\vskip .1in

 In this long example, we will look at the constant sheaf on a space which is not a local complete intersection, and consider a function with a non-isolated critical locus (in any sense of the term ``critical locus''). This example is complicated enough to be interesting, and yet simple enough that we can calculate not only the L\^eVo modules, but also the stalk cohomology of the perverse cohomology of the vanishing cycles. Thus, we can ``check'' our L\^eVo module calculations. 
 
In harder examples, the calculation of the L\^eVo modules would proceed in a similar fashion; though the intersection cycles may be significantly more difficult to compute. Nonetheless, such calculations can be carried out by computer algebra packages. While this falls short of calculating the actual stalk cohomology of the perverse cohomology of the vanishing cycles, it should be viewed as a reasonable --  and effectively calculable -- approximation.

\vskip .2in

Use $(u, x,y,z)$ as
coordinates for
$\Cal U:=\Bbb C^4$, and let $X:=V(u,x)\cup V(y,z)$. Thus, $X$ is the simplest non-local complete intersection. The $\Fdot$ that we will use is $\Bbb Z^\bullet_{X}$. 

Whenever we suppress the reference to the complex of sheaves below, it is assumed that we are using constant $\Bbb Z$ coefficients, and -- in this case -- we will write ordinary cohomology in place of hypercohomology. We continue to use $(w_0, w_1, w_2, w_3)$ as the cotangent coordinates.

\vskip .1in

The strata of $X$ are $S_0:= \{\bold 0\}$, $S_1:= V(u,x)-\{\bold 0\}$, and $S_2:=V(y,z)-\{\bold 0\}$. Obviously, the dimensions of the strata are $d_0=0$, $d_1=2$, and $d_2=2$, respectively. The normal slices to all strata are contractible, while the
corresponding complex links are given by: $\Bbb L_0=$  two complex disks (sets of complex dimension one), $\Bbb L_1=\emptyset$, and $\Bbb L_2=\emptyset$.  Thus, 

\vskip .1in

$H^{k-0}(\Bbb N_0, \Bbb L_0)=0$ unless $k=1$, and  $H^{1-0}(\Bbb N_0, \Bbb L_0)\cong \Bbb Z$;

\vskip .1in

$H^{k-2}(\Bbb N_1, \Bbb L_1)=0$ unless $k=2$, and  $H^{2-2}(\Bbb N_1, \Bbb L_1)\cong \Bbb Z$; and

\vskip .1in

$H^{k-2}(\Bbb N_2, \Bbb L_2)=0$ unless $k=2$, and  $H^{2-2}(\Bbb N_2, \Bbb L_2)\cong \Bbb Z$.

\vskip .1in

\noindent Therefore, $\gecc^k(\Bbb Z^\bullet_X)$ is zero unless $k=1$ or $k=2$, and we find that 
$$
\gecc^1(\Bbb Z^\bullet_X) = \Bbb Z\,[T^*_{\{\bold 0\}}\Cal U] = \Bbb Z\,\left[V(u,x,y,z)\right],$$
$$
\gecc^2(\Bbb Z^\bullet_X) = \Bbb Z\,\left[\overline{T^*_{{}_{S_1}}\Cal
U}\right]+\Bbb Z\,\left[\overline{T^*_{{}_{S_2}}\Cal U}\right] = \Bbb Z\,\left[V(u,x,w_2,w_3)\right]
+\Bbb Z\,\left[V(y,z,w_0,w_1)\right].
$$

\vskip .3in

Now, let  $\tilde f:=
(u^\alpha+x^\beta)^\tau+y^\gamma+z^\delta$,
where $\alpha, \beta,\gamma,\delta,\tau\geqslant 2$. 

\vskip .2in

We calculate as described in Corollary 6.8. We will have to perform two separate such calculations: one where $k=1$ and one where $k=2$. We will follow the notation from Corollary 6.8, except that we will include a superscript on the left to indicate the value of $k$.

\vskip .1in

First, we find that
$$
\imdf = V\big(w_0-\tau\alpha(u^\alpha+x^\beta)^{\tau -1}u^{\alpha-1},  \ 
w_1-\tau\beta(u^\alpha+x^\beta)^{\tau -1}x^{\beta-1},  \ 
w_2- \gamma y^{\gamma-1}, \ w_3-\delta z^{\delta-1}\big)
$$
and so
$$
\imdf\cap|\gecc^1(\Bbb Z^\bullet_X)| = V(u, x, y, z, w_0, w_1, w_2, w_3)
$$
and
$$
\imdf\cap|\gecc^2(\Bbb Z^\bullet_X)| = $$
$$V(u, x, w_2, w_3, w_0, w_1, y, z)\cup V(y, z, w_0, w_1, u^\alpha+x^\beta, w_2, w_3) = V(u^\alpha+x^\beta, y, z, w_0, w_1, w_2, w_3).
$$
As $\dsize\bigcup_{v\in\Bbb C}\supp(\phi_{f-v}[-1]\Bbb Z^\bullet_X) = \eta(|\gecc^\bullet(\Fdot)|\cap\imdf)$, we see that $\phi_{f-v}[-1]\Bbb Z^\bullet_X$ is supported only where $v=0$, and  $\supp(\phi_{f}[-1]\Bbb Z^\bullet_X)  = V(u^\alpha+x^\beta, y, z)$. Thus, $\dim_\bold 0\supp(\phi_f[-1]\Bbb Z^\bullet_X)=1$ and, as mentioned in Remark 6.9, we may calculate the L\^eVo modules as in Corollary 6.8, without worrying ahead of time about whether the coordinates are generic enough; if the intersections in the calculations are all proper, then the coordinates are generic enough.

\vskip .2in

We proceed with the calculation as described in 6.8.

\vskip .1in

Let $${}^1\Pi^4 := \gecc^1(\Bbb Z^\bullet_X) = \Bbb Z\,\left[V(u,x,y,z)\right].$$
 Then,
$$
{}^1\Pi^4\odot V\left(w_3-\frac{\partial\tilde f}{\partial z}\right) \ = \ \Bbb Z\,\left[V(u,x,y,z)\right]\odot [V(w_3-\delta z^{\delta-1})] = \Bbb Z\, [V(u, x, y, z, w_3)],
$$
which has no components contained in $\imdf\cap|\gecc^1(\Bbb Z^\bullet_X)| = V(u, x, y, z, w_0, w_1, w_2, w_3)$. Thus,
$$
{}^1\Pi^3 = \Bbb Z\, [V(u, x, y, z, w_3)],
$$
and we continue, to find
$$
{}^1\Pi^3\odot V\left(w_2-\frac{\partial\tilde f}{\partial y}\right) \ = \ \Bbb Z\, [V(u, x, y, z, w_3)]\odot [V(w_2- \gamma y^{\gamma-1})] \ = \ \Bbb Z\, [V(u, x, y, z, w_3, w_2)] = {}^1\Pi^2.
$$
$$
{}^1\Pi^2\odot V\left(w_1-\frac{\partial\tilde f}{\partial x}\right) \ = \ \Bbb Z\, [V(u, x, y, z, w_3, w_2)]\odot V(w_1-\tau\beta(u^\alpha+x^\beta)^{\tau -1}x^{\beta-1}) =
$$
$$
\Bbb Z\, [V(u, x, y, z, w_3, w_2, w_1)] = {}^1\Pi^1.
$$
$$
{}^1\Pi^1\odot V\left(w_0-\frac{\partial\tilde f}{\partial u}\right)  \ = \ 
\Bbb Z\, [V(u, x, y, z, w_3, w_2, w_1)]\odot V\left(w_0-\tau\alpha(u^\alpha+x^\beta)^{\tau -1}u^{\alpha-1}\right)=
$$
$$
\Bbb Z\, [V(u, x, y, z, w_3, w_2, w_1, w_0)] ={}^1\Delta^0,
$$
where this last equality follows from the fact that $\imdf\cap|\gecc^1(\Bbb Z^\bullet_X)| = V(u, x, y, z, w_0, w_1, w_2, w_3)$. Note that the  above calculations show that ${}^1\Delta^j=0$ if $j\neq 0$.

\vskip .2in

The calculations for $k=2$ are more interesting. Let 
$${}^2\Pi^4 := \gecc^2(\Bbb Z^\bullet_X) =  \Bbb Z\,\left[V(u,x,w_2,w_3)\right]
+\Bbb Z\,\left[V(y,z,w_0,w_1)\right].
$$
$$
{}^2\Pi^4\odot V\left(w_3-\frac{\partial\tilde f}{\partial z}\right) \ = \ 
\Big(\Bbb Z\,\left[V(u,x,w_2,w_3)\right]
+\Bbb Z\,\left[V(y,z,w_0,w_1)\right]\Big)\odot [V(w_3-\delta z^{\delta-1})] = 
$$
$$
(\Bbb Z^{\delta-1})\,\left[V(u,x,w_2,w_3, z)\right] + \Bbb Z\,\left[V(y,z,w_0,w_1, w_3)\right].
$$
As neither of these summands is contained in $\imdf\cap|\gecc^2(\Bbb Z^\bullet_X)|$, we find that $${}^2\Pi^3= (\Bbb Z^{\delta-1})\,\left[V(u,x,w_2,w_3, z)\right] + \Bbb Z\,\left[V(y,z,w_0,w_1, w_3)\right]$$ and we continue.
$$
{}^2\Pi^3\odot V\left(w_2-\frac{\partial\tilde f}{\partial y}\right) \ = \ 
\Big((\Bbb Z^{\delta-1})\,\left[V(u,x,w_2,w_3, z)\right] + \Bbb Z\,\left[V(y,z,w_0,w_1, w_3)\right]\Big)\odot [V(w_2- \gamma y^{\gamma-1})] =
$$
$$
(\Bbb Z^{(\delta-1)(\gamma-1)})\,\left[V(u, x, w_2, w_3, z, y)\right] +\Bbb Z\,\left[V(y, z, w_0, w_1, w_3, w_2)\right] = {}^2\Pi^2.
$$
$$
{}^2\Pi^2\odot V\left(w_1-\frac{\partial\tilde f}{\partial x}\right) =
$$
$$
\Big((\Bbb Z^{(\delta-1)(\gamma-1)})\,\left[V(u, x, w_2, w_3, z, y)\right] +\Bbb Z\,\left[V(y, z, w_0, w_1, w_3, w_2)\right]\Big)\odot V(w_1-\tau\beta(u^\alpha+x^\beta)^{\tau -1}x^{\beta-1}) =
$$
$$
(\Bbb Z^{(\delta-1)(\gamma-1)})\,\left[V(u, x, w_2, w_3, z, y, w_1)\right] + \Bbb Z\,\left[V(y, z, w_0,w_1, w_3, w_2, (u^\alpha+x^\beta)^{\tau -1}x^{\beta-1})\right] =
$$
$$
(\Bbb Z^{(\delta-1)(\gamma-1)})\,\left[V(u,x,w_2,w_3, z, y, w_1)\right] + ( \Bbb Z^{\beta-1})\,\left[V(y, z, w_0, w_1, w_3, w_2, x)\right] +
$$
$$
 ( \Bbb Z^{\tau-1})\,\left[V(y, z, w_0, w_1, w_3, w_2, u^\alpha+x^\beta)\right],
$$
where this last summand is contained in $\imdf$, but the two earlier ones were not. Thus,
$$
{}^2\Delta^1= ( \Bbb Z^{\tau-1})\,\left[V(y, z, w_0, w_1, w_3, w_2, u^\alpha+x^\beta)\right]
$$
and
$$
{}^2\Pi^1= (\Bbb Z^{(\delta-1)(\gamma-1)})\,\left[V(u, x, w_2, w_3, z, y, w_1)\right] + ( \Bbb Z^{\beta-1})\,\left[V(y, z, w_0, w_1, w_3, w_2, x)\right] .$$
Finally, we find
$$
{}^2\Pi^1\odot V\left(w_0-\frac{\partial\tilde f}{\partial u}\right) =
$$
$$
\Big((\Bbb Z^{(\delta-1)(\gamma-1)})\,\left[V(u, x, w_2, w_3, z, y, w_1)\right] + ( \Bbb Z^{\beta-1})\,\left[V(y, z, w_0, w_1, w_3, w_2, x)\right]\Big)\odot 
$$
$$V(w_0-\tau\alpha(u^\alpha+x^\beta)^{\tau -1}u^{\alpha-1}) =
$$
$$
(\Bbb Z^{(\delta-1)(\gamma-1)})\,\left[V(u, x, w_2, w_3, z, y, w_1, w_0)\right] +
( \Bbb Z^{(\beta-1)(\alpha-1)})\,\left[V(y, z, w_0, w_1, w_3, w_2, x, u)\right]+$$
$$
( \Bbb Z^{(\beta-1)(\tau-1)\alpha})\,\left[V(y, z, w_0, w_1, w_3, w_2, x, u)\right] = 
$$
$$
( \Bbb Z^{(\delta-1)(\gamma-1)+(\beta-1)(\alpha\tau-1)})\,\left[V(u, x, y, z, w_0,w_1, w_2, w_3)\right] ={}^2\Delta^0.
$$

Projecting down into $\Cal U$ (and suppressing the reference to the coordinates and to $\Bbb Z^\bullet_X$), we find that the only non-zero L\^eVo cycles are:
$$
{}^1\Lambda^0_f = \Bbb Z\,[\{\bold 0\}],
$$
$$
{}^2\Lambda^1_f= ( \Bbb Z^{\tau-1})\,\left[V(u^\alpha+x^\beta, y, z)\right],
$$
and
$$
{}^2\Lambda^0_f= ( \Bbb Z^{(\delta-1)(\gamma-1)+(\beta-1)(\alpha\tau-1)})\,[\{\bold 0\}].$$

The corresponding L\^eVo modules at the origin are 
$$
{}^1\lambda^0_f(\bold 0)=\Bbb Z
$$
$$
{}^2\lambda^1_f(\bold 0)=\Big(V(u)\odot {}^2\Lambda^1_f\Big)_{\{\bold 0\}}= \Big(V(u)\odot( \Bbb Z^{\tau-1})\,\left[V(u^\alpha+x^\beta, y, z)\right]\Big)_{\{\bold 0\}}=\Bbb Z^{\beta(\tau-1)},
$$
and
$$
{}^2\lambda^0_f(\bold 0)=  \Bbb Z^{(\delta-1)(\gamma-1)+(\beta-1)(\alpha\tau-1)}.$$

\vskip .3in

We will now calculate the stalk cohomology of the perverse cohomology of the vanishing cycles, and then compare that data with the L\^eVo module data. This comparison will use the Zawatsky complex (5.19) for the vanishing cycles, which was the next-to-last statement in Theorem 6.4.

\vskip .1in

For every locally closed subset $Z\subseteq X$, we let $(\Bbb Z^\bullet _X)_Z$ denote the restriction-extension of $\Bbb Z^\bullet _X$ on $Z$, i.e., the complex obtained by first restricting to $Z$ and then by extending by zero over all of $X$.

 Let $X_1:= V(u,x)$ and $X_2:=V(y, z)$. Then, there is a canonical distinguished triangle
$$
\Bbb Z^\bullet _X\rightarrow (\Bbb Z^\bullet _X)_{X_1}\oplus(\Bbb Z^\bullet _X)_{X_2}\rightarrow(\Bbb Z^\bullet _X)_{\{\bold 0\}}@>[1]>>\Bbb Z^\bullet _X.
$$
Note that $(\Bbb Z^\bullet _X)_{X_1}[2]$, $(\Bbb Z^\bullet _X)_{X_1}[2]$, and  $(\Bbb Z^\bullet _X)_{\{\bold 0\}}$ are all perverse, and so the only non-zero portion of the long exact sequence on perverse cohomology resulting from the above distinguished triangle becomes
$$
 0\rightarrow {}^\mu\hskip -.02in H^0(\Bbb Z^\bullet _X)\rightarrow 0\rightarrow (\Bbb Z^\bullet _X)_{\{\bold 0\}}\rightarrow  {}^\mu\hskip -.02in H^1(\Bbb Z^\bullet _X)\rightarrow 0\rightarrow 0\rightarrow  {}^\mu\hskip -.02in H^2(\Bbb Z^\bullet _X)\rightarrow (\Bbb Z^\bullet _X)_{X_1}\oplus(\Bbb Z^\bullet _X)_{X_2}\rightarrow 0.
$$
Therefore, ${}^\mu\hskip -.02in H^k(\Bbb Z^\bullet _X)$ is unequal to zero only when $k=1$ or $k=2$, and
${}^\mu\hskip -.02in H^1(\Bbb Z^\bullet _X)\cong (\Bbb Z^\bullet _X)_{\{\bold 0\}}$ and ${}^\mu\hskip -.02in H^2(\Bbb Z^\bullet _X)\cong  (\Bbb Z^\bullet _X)_{X_1}\oplus(\Bbb Z^\bullet _X)_{X_2}$.

It follows that $$\phi_f[-1]\big({}^\mu\hskip -.02in H^1(\Bbb Z^\bullet _X)\big)\cong {}^\mu\hskip -.02in H^1(\phi_f[-1]\Bbb Z^\bullet _X)\cong (\Bbb Z^\bullet _X)_{\{\bold 0\}}$$
and
 $$\phi_f[-1]\big({}^\mu\hskip -.02in H^2(\Bbb Z^\bullet _X)\big)\cong {}^\mu\hskip -.02in H^2(\phi_f[-1]\Bbb Z^\bullet _X)\cong \phi_f[-1]\big((\Bbb Z^\bullet _X)_{X_1}\big)\oplus\phi_f[-1]\big((\Bbb Z^\bullet _X)_{X_2}\big).
$$
If we let $f_1:= f_{|_{X_1}}$, $f_2:= f_{|_{X_2}}$, and let $j_1: X_1\hookrightarrow X$ and $j_2: X_2 \hookrightarrow X$ denote the inclusions, then we can rewrite the last isomorphisms above as
$$\phi_f[-1]\big({}^\mu\hskip -.02in H^2(\Bbb Z^\bullet _X)\big)\cong {}^\mu\hskip -.02in H^2(\phi_f[-1]\Bbb Z^\bullet _X)\cong {j_1}_!\big(\phi_{f_1}[-1]\Bbb Z^\bullet _{X_1}\big)\oplus{j_2}_!\big(\phi_{f_2}[-1]\Bbb Z^\bullet _{X_2}\big).
$$

\vskip .1in

Now, $F_1$, the Milnor fibre of $f_1$ at the origin, is the Milnor fibre of $y^\gamma+z^\delta$ in the $yz$-plane, and $F_2$, the Milnor fibre of $f_2$ at the origin, is the Milnor fibre of $(u^\alpha+x^\beta)^\tau$ in the $ux$-plane. Thus, $F_1$ is homotopy-equivalent to a bouquet of $(\gamma-1)(\delta-1)$ circles, and $F_2$
is homotopy-equivalent to the disjoint union of $\tau$ copies of a bouquet of $(\alpha-1)(\beta-1)$ circles. 

\vskip .1in

We conclude from all of this that the stalk cohomology module $H^j\big({}^\mu\hskip -.02in H^k(\phi_f[-1]\Bbb Z^\bullet _X)\big)_\bold 0$ is zero except for the following:
$$
H^0\big({}^\mu\hskip -.02in H^1(\phi_f[-1]\Bbb Z^\bullet _X)\big)_\bold 0\cong\Bbb Z,
$$
$$
H^{-1}\big({}^\mu\hskip -.02in H^2(\phi_f[-1]\Bbb Z^\bullet _X)\big)_\bold 0\cong \Bbb Z^{\tau-1},
$$
and
$$
H^0\big({}^\mu\hskip -.02in H^2(\phi_f[-1]\Bbb Z^\bullet _X)\big)_\bold 0\cong \Bbb Z^{(\gamma-1)(\delta-1)+\tau(\alpha-1)(\beta-1)}.
$$

\vskip .2in

Now, the Zawatsky complex, together with our earlier calculation of the L\^eVo modules, also tells us that only the above stalk cohomology modules could possibly be non-zero. Furthermore, the degree $1$ Zawatsky complex implies that we had to have $H^0\big({}^\mu\hskip -.02in H^1(\phi_f[-1]\Bbb Z^\bullet _X)\big)_\bold 0\cong {}^1\lambda^0_f(\bold 0)\cong \Bbb Z$. 

The degree $2$ Zawatsky complex tells us that there is a homomorphism  $$\omega: \Bbb Z^{\beta(\tau-1)}\cong {}^2\lambda^1_f(\bold 0) \ \longrightarrow \  {}^2\lambda^0_f(\bold 0)\cong  \Bbb Z^{(\delta-1)(\gamma-1)+(\beta-1)(\alpha\tau-1)}$$
such that $\operatorname{ker}\omega \cong H^{-1}\big({}^\mu\hskip -.02in H^2(\phi_f[-1]\Bbb Z^\bullet _X)\big)_\bold 0$ and $\operatorname{coker}\omega\cong H^0\big({}^\mu\hskip -.02in H^2(\phi_f[-1]\Bbb Z^\bullet _X)\big)_\bold 0$.

\vskip .1in

This is consistent with our perverse cohomology calculations since:

\vskip .1in

\noindent $\bullet$\hskip .2in $H^{-1}\big({}^\mu\hskip -.02in H^2(\phi_f[-1]\Bbb Z^\bullet _X)\big)_\bold 0$ is free, of rank no more than $\beta(\tau-1)$;

\vskip .1in

\noindent $\bullet$\hskip .2in $H^{0}\big({}^\mu\hskip -.02in H^2(\phi_f[-1]\Bbb Z^\bullet _X)\big)_\bold 0$ has rank no more than $(\delta-1)(\gamma-1)+(\beta-1)(\alpha\tau-1)$; and

\vskip .1in

\noindent $\bullet$\hskip .2in the alternating sums agree, i.e.,
$$
\Big((\delta-1)(\gamma-1)+(\beta-1)(\alpha\tau-1)\Big) - \Big(\beta(\tau-1)\Big) \ = \ \Big((\gamma-1)(\delta-1)+\tau(\alpha-1)(\beta-1)\Big)-\Big(\tau-1\Big).
$$

\vskip .1in

Finally, the reader should check that either method yields the same Euler characteristic for the reduced cohomology of the Milnor fibre of $f$ at the origin; namely, 
$$
\widetilde\chi\big(F_{f, \bold 0}\big) =-\alpha\beta\tau+\beta\tau+\alpha\tau-\gamma\delta+\gamma+\delta-1.
$$

\vskip .3in

\noindent\S8. {\bf Remarks and Future Directions }

\vskip .2in

What we have shown in this paper is that the characteristic polar modules provide an enriched form of the classical absolute polar multiplicities, and the L\^e-Vogel modules provide an enriched form of the L\^e numbers. Moreover, these enriched pieces of data encode a great deal of data associated to complexes of sheaves, and yet these enriched devices remain calculable in an algebraic/intersection-theoretic manner.

\vskip .1in

Originally, this paper was to be entitled ``Enriched Cycles and Equisingularity''. However, as the list of fundamental definitions and results grew, it became clear that most of the equisingularity results would have to go into another paper. We say ``most'' because we view Theorem 6.5, and even Corollary 3.5, as equisingularity results.

What we will prove in Enriched Cycles and Equisingularity is that the constancy of the L\^e-Vogel modules throughout a family implies that the $a_f$ condition holds, and that the stalk cohomology of the vanishing cycles is constant.

As fundamental tools in proving these results, we will prove an upper-semicontinuity result and general L\^e-Iomdine-Vogel formulas. The notion of upper-semicontinuity makes sense for enriched cycles due to the existence of the partial ordering on enriched cycles given in Section 2. The L\^e-Iomdine-Vogel formulas are a generalization of the L\^e-Iomdine formulas of [{\bf M3}], and use the general Vogel cycle results that we developed in Part I of [{\bf M7}]; these formulas inductively reduce a number of problems to the case of isolated critical points.

\vskip .3in

In the future, we will also investigate the extent to which one can encode monodromy results, and prove monodromy theorems, via enriched cycles;  we wish to capture many aspects of the monodromy on $\psi_f[-1]\Fdot$ and $\phi_f[-1]\Fdot$.

It may appear that we have not given ourselves enough structure to discuss the monodromies -- for we actually defined two enriched cycles to be equal if the component modules were isomorphic. In other words, an enriched cycle is actually an equivalence class. This was necessary since there are no canonical choices for defining the Morse data to strata.

However, we can also define the monodromies as equivalence classes. If we fix one representative for an enriched cycle, the monodromy is determined by an automorphism of each component module. If we switch to a different representative of the enriched cycle, an equivalent family of automorphisms is one obtained by pulling-back the original family of automorphisms via the isomorphisms of the component modules, i.e., the new automorphisms are conjugate to the original ones.

It will be interesting to see to what extent we can recover portions of such deep results as the monodromy theorem and the decomposition theorem through enriched cycle techniques.

\vskip .3in

Finally, it would greatly enhance the algorithmic aspect of our work if we could prove Proposition 6.11 regardless of the dimension of $\Sigma_{{}_{\Fdot}}f$. On the other hand, it would also be interesting to find a counterexample for higher-dimensional critical loci.

\enddocument